\documentclass[11pt]{article}

\begin{document}



\newcommand{\newc}{\newcommand}


\renewcommand{\theequation}{\thesection.\arabic{equation}}
\newc{\eqnoset}{\setcounter{equation}{0}}
\newcommand{\myref}[2]{#1~\ref{#2}}

\newcommand{\mref}[1]{(\ref{#1})}
\newcommand{\reflemm}[1]{Lemma~\ref{#1}}
\newcommand{\refrem}[1]{Remark~\ref{#1}}
\newcommand{\reftheo}[1]{Theorem~\ref{#1}}
\newcommand{\refdef}[1]{Definition~\ref{#1}}
\newcommand{\refcoro}[1]{Corollary~\ref{#1}}
\newcommand{\refprop}[1]{Proposition~\ref{#1}}
\newcommand{\refsec}[1]{Section~\ref{#1}}
\newcommand{\refchap}[1]{Chapter~\ref{#1}}

\newcommand{\beq}{\begin{equation}}
\newcommand{\eeq}{\end{equation}}
\newcommand{\beqno}[1]{\begin{equation}\label{#1}}

\newcommand{\barr}{\begin{array}}
\newcommand{\earr}{\end{array}}

\newc{\bearr}{\begin{eqnarray*}}
\newc{\eearr}{\end{eqnarray*}}

\newc{\bearrno}[1]{\begin{eqnarray}\label{#1}}
\newc{\eearrno}{\end{eqnarray}}

\newc{\non}{\nonumber}
\newc{\nol}{\nonumber\nl}

\newcommand{\bdes}{\begin{description}}
\newcommand{\edes}{\end{description}}
\newc{\benu}{\begin{enumerate}}
\newc{\eenu}{\end{enumerate}}
\newc{\btab}{\begin{tabular}}
\newc{\etab}{\end{tabular}}



\newtheorem{theorem}{Theorem}[section]
\newtheorem{defi}[theorem]{Definition}
\newtheorem{lemma}[theorem]{Lemma}
\newtheorem{rem}[theorem]{Remark}
\newtheorem{exam}[theorem]{Example}
\newtheorem{propo}[theorem]{Proposition}
\newtheorem{corol}[theorem]{Corollary}

\renewcommand{\thelemma}{\thesection.\arabic{lemma}}

\newcommand{\btheo}[1]{\begin{theorem}\label{#1}}
\newc{\brem}[1]{\begin{rem}\label{#1}\em}
\newc{\bexam}[1]{\begin{exam}\label{#1}\em}
\newc{\bdefi}[1]{\begin{defi}\label{#1}}
\newcommand{\blemm}[1]{\begin{lemma}\label{#1}}
\newcommand{\bprop}[1]{\begin{propo}\label{#1}}
\newcommand{\bcoro}[1]{\begin{corol}\label{#1}}
\newcommand{\etheo}{\end{theorem}}
\newcommand{\elemm}{\end{lemma}}
\newcommand{\eprop}{\end{propo}}
\newcommand{\ecoro}{\end{corol}}
\newc{\erem}{\end{rem}}
\newc{\eexam}{\end{exam}}
\newc{\edefi}{\end{defi}}

\newc{\rmk}[1]{{\bf REMARK #1: }}
\newc{\DN}[1]{{\bf DEFINITION #1: }}

\newcommand{\bproof}{{\bf Proof:~~}}
\newc{\eproof}{{\vrule height8pt width5pt depth0pt}\vspace{3mm}}


\newcommand{\rarr}{\rightarrow}
\newcommand{\Rarr}{\Rightarrow}
\newcommand{\tru}{\backslash}
\newc{\bfrac}[2]{\dspl{\frac{#1}{#2}}}


\newc{\nl}{\vspace{2mm}\\}
\newc{\nid}{\noindent}


\newcommand{\oneon}[1]{\frac{1}{#1}}
\newcommand{\dspl}{\displaystyle}
\newc{\grad}{\nabla}
\newc{\Div}{\mbox{div}}
\newc{\pdt}[1]{\dspl{\frac{\partial{#1}}{\partial t}}}
\newc{\pdn}[1]{\dspl{\frac{\partial{#1}}{\partial \nu}}}
\newc{\pdNi}[1]{\dspl{\frac{\partial{#1}}{\partial \mathcal{N}_i}}}
\newc{\pD}[2]{\dspl{\frac{\partial{#1}}{\partial #2}}}
\newc{\dt}{\dspl{\frac{d}{dt}}}
\newc{\bdry}[1]{\mbox{$\partial #1$}}
\newc{\sgn}{\mbox{sign}}

\newc{\Hess}[1]{\frac{\partial^2 #1}{\pdh z_i \pdh z_j}}
\newc{\hess}[1]{\partial^2 #1/\pdh z_i \pdh z_j}


\newcommand{\Coone}[1]{\mbox{$C^{1}_{0}(#1)$}}
\newcommand{\lspac}[2]{\mbox{$L^{#1}(#2)$}}
\newc{\hspac}[2]{\mbox{$C^{0,#1}(#2)$}}
\newc{\Hspac}[2]{\mbox{$C^{1,#1}(#2)$}}
\newc{\Hosp}{\mbox{$H^{1}_{0}$}}
\newcommand{\Lsp}[1]{\mbox{$L^{#1}(\Og)$}}
\newc{\hsp}{\Hosp(\Og)}


\newc{\ag}{\alpha}
\newc{\bg}{\beta}
\newc{\cg}{\gamma}\newc{\Cg}{\Gamma}
\newc{\dg}{\delta}\newc{\Dg}{\Delta}
\newc{\eg}{\varepsilon}
\newc{\zg}{\zeta}
\newc{\thg}{\theta}
\newc{\llg}{\lambda}\newc{\LLg}{\Lambda}
\newc{\kg}{\kappa}
\newc{\rg}{\rho}
\newc{\sg}{\sigma}\newc{\Sg}{\Sigma}
\newc{\tg}{\tau}
\newc{\fg}{\phi}\newc{\Fg}{\Phi}
\newc{\vfg}{\varphi}
\newc{\og}{\omega}\newc{\Og}{\Omega}
\newc{\pdh}{\partial}


\newc{\ii}[1]{\int_{#1}}
\newc{\iidx}[2]{{\dspl\int_{#1}~#2~dx}}
\newc{\bii}[1]{{\dspl \ii{#1} }}
\newc{\biii}[2]{{\dspl \iii{#1}{#2} }}
\newc{\su}[2]{\sum_{#1}^{#2}}
\newc{\bsu}[2]{{\dspl \su{#1}{#2} }}

\newcommand{\iiomdx}[1]{{\dspl\int_{\Og}~ #1 ~dx}}
\newc{\biiom}[1]{{\dspl\int_{\bdrom}~ #1 ~d\sg}}
\newc{\io}[1]{{\dspl\int_{\Og}~ #1 ~dx}}
\newc{\bio}[1]{{\dspl\int_{\bdrom}~ #1 ~d\sg}}
\newc{\bsir}{\bsu{i=1}{r}}
\newc{\bsim}{\bsu{i=1}{m}}

\newc{\iibr}[2]{\iidx{\bprw{#1}}{#2}}
\newc{\Intbr}[1]{\iibr{R}{#1}}
\newc{\intbr}[1]{\iibr{\rg}{#1}}
\newc{\intt}[3]{\int_{#1}^{#2}\int_\Og~#3~dxdt}

\newc{\itQ}[2]{\dspl{\int\hspace{-2.5mm}\int_{#1}~#2~dz}}
\newc{\mitQ}[2]{\dspl{\rule[1mm]{4mm}{.3mm}\hspace{-5.3mm}\int\hspace{-2.5mm}\int_{#1}~#2~dz}}
\newc{\mitQQ}[3]{\dspl{\rule[1mm]{4mm}{.3mm}\hspace{-5.3mm}\int\hspace{-2.5mm}\int_{#1}~#2~#3}}

\newc{\mitx}[2]{\dspl{\rule[1mm]{3mm}{.3mm}\hspace{-4mm}\int_{#1}~#2~dx}}

\newc{\mitQq}[2]{\dspl{\rule[1mm]{4mm}{.3mm}\hspace{-5.3mm}\int\hspace{-2.5mm}\int_{#1}~#2~d\bar{z}}}
\newc{\itQq}[2]{\dspl{\int\hspace{-2.5mm}\int_{#1}~#2~d\bar{z}}}

\newc{\pder}[2]{\dspl{\frac{\partial #1}{\partial #2}}}


\newc{\ui}{u_{i}}
\newcommand{\upl}{u^{+}}
\newcommand{\umn}{u^{-}}
\newcommand{\un}{\{ u_{n}\}}

\newcommand{\uo}{u_{0}}
\newc{\voi}{v_{i}^{0}}
\newc{\uoi}{u_{i}^{0}}
\newc{\vu}{\vec{u}}

\newc{\xo}{x_{0}}
\newc{\Br}{B_{R}}
\newc{\Bro}{\Br (\xo)}
\newc{\bdrom}{\bdry{\Og}}
\newc{\ogr}[1]{\Og_{#1}}
\newc{\Bxo}{B_{x_0}}

\newc{\inP}[2]{\|#1(\bullet,t)\|_#2\in\cP}
\newc{\cO}{{\mathcal O}}
\newc{\inO}[2]{\|#1(\bullet,t)\|_#2\in\cO}

\newc{\newl}{\\ &&}

\newc{\bilhom}{\mbox{Bil}(\mbox{Hom}(\RR^{nm},\RR^{nm}))}
\newc{\VV}[1]{{V(Q_{#1})}}

\newc{\ccA}{{\mathcal A}}
\newc{\ccB}{{\mathcal B}}
\newc{\ccC}{{\mathcal C}}
\newc{\ccD}{{\mathcal D}}
\newc{\ccE}{{\mathcal E}}
\newc{\ccH}{\mathcal{H}}
\newc{\ccG}{\mathcal{G}}
\newc{\ccF}{\mathcal{F}}
\newc{\ccL}{\mathcal{L}}
\newc{\ccK}{\mathcal{K}}
\newc{\ccI}{{\mathcal I}}
\newc{\ccJ}{{\mathcal J}}
\newc{\ccP}{{\mathcal P}}
\newc{\ccQ}{{\mathcal Q}}
\newc{\ccR}{{\mathcal R}}
\newc{\ccS}{{\mathcal S}}
\newc{\ccT}{{\mathcal T}}
\newc{\ccX}{{\mathcal X}}
\newc{\ccY}{{\mathcal Y}}
\newc{\ccZ}{{\mathcal Z}}

\newc{\bb}[1]{{\mathbf #1}}
\newc{\bbA}{{\mathbf A}}
\newc{\myprod}[1]{\langle #1 \rangle}
\newc{\mypar}[1]{\left( #1 \right)}


\newc{\lspn}[2]{\mbox{$\| #1\|_{\Lsp{#2}}$}}
\newc{\Lpn}[2]{\mbox{$\| #1\|_{#2}$}}
\newc{\Hn}[1]{\mbox{$\| #1\|_{H^1(\Og)}$}}


\newc{\cyl}[1]{\og\times \{#1\}}
\newc{\cyll}{\og\times[0,1]}
\newc{\vx}[1]{v\cdot #1}
\newc{\vtx}[1]{v(t,x)\cdot #1}
\newc{\vn}{\vx{n}}

\newcommand{\RR}{{\rm I\kern -1.6pt{\rm R}}}

\newc{\mX}{{\mathbf{X}}}
\newc{\mP}{{\mathbf{P}}}
\newc{\mB}{{\mathbf{B}}}


\newenvironment{proof}{\noindent\textbf{Proof.}\ }
{\nopagebreak\hbox{ }\hfill$\Box$\bigskip}


\newc{\itQQ}[2]{\dspl{\int_{#1}#2\,dz}}
\newc{\mmitQQ}[2]{\dspl{\rule[1mm]{4mm}{.3mm}\hspace{-4.3mm}\int_{#1}~#2~dz}}
\newc{\MmitQQ}[2]{\dspl{\rule[1mm]{4mm}{.3mm}\hspace{-4.3mm}\int_{#1}~#2~d\mu}}

\newc{\MUmitQQ}[3]{\dspl{\rule[1mm]{4mm}{.3mm}\hspace{-4.3mm}\int_{#1}~#2~d#3}}
\newc{\MUitQQ}[3]{\dspl{\int_{#1}~#2~d#3}}

\vspace*{-.8in}
\begin{center} {\LARGE\em  Existence of Strong and Nontrivial Solutions to Strongly Coupled Elliptic Systems}

 \end{center}

\vspace{.1in}

\begin{center}

{\sc Dung Le}{\footnote {Department of Mathematics, University of
Texas at San
Antonio, One UTSA Circle, San Antonio, TX 78249. {\tt Email: Dung.Le@utsa.edu}
}}

\end{center}

\begin{abstract}
We establish the existence of strong solutions to a class of nonlinear strongly coupled and uniform  elliptic systems consisting of more than two equations. The existence of of nontrivial and non constant solutions (or pattern formations) will also be studied.
\end{abstract}

\vspace{.2in}

\section{Introduction}\label{introsec}\eqnoset

In this paper, we study the existence of {\em strong} solutions and {\em other nontrivial} solutions to the following nonlinear {\em strongly coupled} and {\em  nonregular} but {\em uniform}  elliptic system
\beqno{e10}\left\{\barr{l} -\Div(A(u,Du))=\hat{f}(u,Du) \mbox{ in $\Og$},\\\mbox{$u$ satisfies Dirichlet or Neumann boundary conditions on $\partial \Og$}. \earr\right.\eeq

Here, $\Og$ is a bounded domain with smooth boundary $\partial \Og$ in $\RR^n$, $n\ge2$. A typical point in $\RR^n$ is denoted by $x$. The $k$-order derivatives of a vector valued function $$u(x)=(u_1(x),\ldots,u_m(x)) \quad m\ge2 $$ are denoted by $D^ku$. $A(u,Du)$ is a {\em full} matrix $m\times n$ and $\hat{f}:\RR^m\times\RR^{nm}\to\RR^m$. 
Also, for a vector or matrix valued function $f(u,\zeta)$, $u\in \RR^m$ and $\zeta\in \RR^d$, its partial derivatives will be denoted by $f_u,f_\zeta$. 

Throughout this paper, we always assume the following on the diffusion matrix $A(u,Du)$.

\bdes

\item[A)] $A(u,\zeta)$ is $C^1$ in   $u\in\RR^m$ and $\zeta\in\RR^{mn}$. There are a constant $C_*>0$  and a nonnegative  scalar $C^1$ function $\llg(u)$  such that  for any  $u\in\RR^m$ and $\zeta,\xi\in\RR^{nm}$ 
\beqno{A1} \llg(u)|\zeta|^2 \le \myprod{A_\zeta(u,\zeta)\xi,\xi} \mbox{ and } |A_\zeta(u,\zeta)|\le C_*\llg(u).\eeq 

Moreover, there are positive constants $C,\llg_0$ such that $\llg(u)\ge \llg_0$ and $|A_u(u,\zeta)|\le C|\llg_u(u)||\zeta|$. In addition, $A(u,0)=0$ for all $u\in\RR^m$.
\edes

The first condition in \mref{A1} is to say that the system \mref{e10} is elliptic. If $\llg(u)$ is also bounded from above by a constant for all $u\in\RR^m$, we say that $A$ is regular elliptic. Otherwise, $A$ is uniform elliptic if \mref{A1} holds.

Furthermore, the constant $C_*$ in \mref{A1} concerns the ratio betweeen the largest and smallest eigenvalues of $A_\zeta$. We assume that these  constants are not too far apart in the following sense.
\bdes\item[SG)] (The spectral gap condition)  $ C_*<(n-2)/(n-4) $ if $n>4$. \edes

We note that if SG) is somewhat violated, i.e. $C_*$ is large, then examples of blowing up in finite time can occur for the corresponding parabolic systems (see \cite{sd}). Of course, this condition is void in many applications when we usually have $n\le4$.

By a {\em strong solution} of \mref{e10} we mean a vector valued function $u\in W^{2,p}(\Og,\RR^m)$ for any $p>1$ that solves \mref{e1} a.e. in $\Og$ and
$Du\in C^{\ag}(\Og,,\RR^m)$, $\ag\in(0,1)$. 

As usual,  $W^{k,p}(\Og,\RR^m)$,  where $k$ is an integer and $p\ge1$,  denotes the standard Sobolev spaces whose elements are vector valued functions $u\,:\,\Og\to \RR^m$ with finite norm $$\|u\|_{W^{k,p}(\Og,\RR^m)} = \|u\|_{L^p(\Og,\RR^m)} + \sum_{i=1}^k\|D^{k}u\|_{L^p(\Og,\RR^{kmn})}.$$ Similarly, $C^{k,\ag}(\Og,\RR^m)$ denotes the space of (vector valued) functions $u$ on $\Og$ such that $D^lu$, $l=0,\ldots,k$, are H\"older continuous with exponent $\ag\in(0,1)$. If the range $\RR^m$ is understood from the context we will usually omit it from the above notations.

The system \mref{e10} occurs in many applications concerning steady states of diffusion processes with cross diffusion taken into account, i.e. $A(u,Du)$ is a full matrix (see \cite{JA} and the reference therein). In the last few decades, there are many studies of \mref{e10} under the main assumption that its solutions are bounded. The lack of maximum principles for systems of more than one equations has limited the range of application of those results. Occasionally, works in this direction usually tried to establish $L^\infty$ bounds for solutions via ad hoc techniques and thus imposed restrictive assumptions on the structural conditions of the systems. On the other hand, even if $L^\infty$ boundedness of solutions were known, counterexamples in \cite{JS} showed that this does not suffice to guarantee higher regularity of the solutions.

Our first goal is to establish the existence of a strong solution to the general \mref{e10} when its $L^\infty$ boundedness is not available. Since the system is not variational and comparison principles are generally unvailable,  techniques in variational methods and monotone dynamical systems  are not applicable here. Fixed point index theories will then be more appropriate to study the existence of solutions to \mref{e10}. However, it is well known that the main ingredients of this approach are: 1) to define compact map $T$,  whose fixed points are solutions to \mref{e10}, on some appropriate Banach space $\mX$; 2) to show that the Leray Schauder fixed point index of $T$ is nonzero. The second part requires some uniform estimates of the fixed points of $T$ and regularity properties of solutions to \mref{e10}. Those are the fundamental and most technical problems in the theory of partial differential equations. In this work, we will show that the crucial regularity property can be obtained if we know a~priori that the solutions are VMO or BMO in the case of large self diffusion. We will show that the result applies to the generalized SKT systems (\cite{SKT}) when the dimension of $\Og$ is less than 5.
 
To this end, and throughout this work, we will impose the following structural conditions on the reaction term $\hat{f}$.

\bdes \item[F)] There exists a constant $C$  such that for any $C^1$ functions $u:\Og\to\RR^m$ and $p:\Og\to\RR^{mn}$ satisfying $|Du|\le |p|$ the following holds.  \beqno{FUDU11a}|D\hat{f}(u,p)| \le C[\llg(u)|Dp| + |\llg_u(u)||Du||p|+\llg(u)||p|].\eeq
\edes

 In addition, we can see that F) typically holds if there exist a constant $C$ and a function $f(u)$ which is $C^1$ in $u$  such that \beqno{FUDU11}|\hat{f}(u,p)| \le C\llg(u)|p| + f(u),\quad |f_u(u)| \le C\llg(u).\eeq

For simplicity, we are assuming a linear growth in $p$ on $\hat{f}$. In fact, the main existence results in this work allow the following nonlinear growth for $\hat{f}$ $$|\hat{f}(u,p)| \le C\llg(u)|p|^\ag + f(u)  \quad \mbox{for some $\ag\in[1,2)$},$$ and such that $$|D\hat{f}(u,p)| \le C[\llg(u)|p|^{\ag-1}|Dp| + |\llg_u(u)||Du||p|^\ag+\llg(u)|p|^\ag].$$ The proof is similar with minor modifications (see \refrem{fdurem1}).

 To establish the existence of a strong solution, we embed \mref{e10} in a suitable family of systems with $\sg\in[0,1]$
 \beqno{famsysintro}\left\{\barr{l} -\Div(\hat{A}_\sg(U,DU))=\hat{F}_\sg(U,DU) \mbox{ in $\Og$,}\\
 \mbox{Homogeneous Dirichlet or Neumann boundary conditions on $\partial \Og$.}
 \earr\right.\eeq The data $\hat{A}_\sg,\hat{F}_\sg$ satisfy A), F) and SG) with the same set of constants. We then consider a family of compact maps $T(\sg,\cdot)$ associated to the above systems and use a homotopy argument to compute the fixed point index of $T$. Again, the key point is to establish some uniform estimates of the fixed points of $T(\sg,\cdot)$ and regularity properties of their fixed points.

 The uniform estimates for H\"older norms and then higher norms of solutions to the above systems come from the crucial and technical \refprop{Dup-uni} in \refsec{estsec} which shows that one needs only a uniform control of the $W^{1,2}(\Og)$ and $VMO(\Og)$ norms of (unbounded) strong solutions to the systems. Roughly speaking, we assume that for any given $\mu_0>0$ there is a positive $R_{\mu_0}$ for which the strong solutions to the systems in \mref{famsysintro} satisfy  \beqno{KeyVMO1} \mathbf{\LLg}^2\sup_{x_0\in\bar{\Og}}\|U\|_{BMO(B_{x_0})}^2 \le \mu_0.\eeq

 The proof of this result relies on a combination of a local weighted Gagliardo-Nirenberg inequality which is proved in our recent work \cite{dleGNans} (see also \cite{letrans}) and a new iteration argument using decay estimates. This technique was used in our work \cite{dleGNans} to establish the global existence of solutions to strongly coupled parabolic systems. The proof for the elliptic case in this paper is somewhat simpler and requires less assumptions but needs some subtle modifications. For the sake of completeness and the convenience of the readers we present the details.
 
 The fact that {\em bounded} weak $VMO$ solutions to {\em regular} elliptic systems are H\"older continuous is now well known (see \cite{Gius}). Here, using a completely different approach, we deal with {\em unbounded strong} $VMO$ solutions and our \refprop{Dup-uni} applies to  {\em nonregular} systems. Eventually, we obtain that the strong solutions to the systems of \mref{famsysintro} are uniformly H\"older continuous. Desired uniform estimates for higher norms of the solutions then follow.
 
  Once this technical result is established, our first main result in \refsec{extsec}, \reftheo{extthm}, then shows that \mref{e10} has a strong solution if the strong solutions of \mref{famsysintro} are uniformly bounded in $W^{1,2}(\Og)$ and $VMO(\Og)$.

 We  present some examples in applications where \reftheo{extthm} can apply. The main theme in these examples is to establish the uniform boundedness of the  of solutions to \mref{e10} in $W^{1,n}(\Og)$, so that the solutions are in $W^{1,2}(\Og)$ and $VMO(\Og)$. In fact, under suitable assumptions, which occur in many mathematical models in biology and ecology, on the structural of \mref{e10} we will show that it is sufficient to control the very weak $L^1$ norms of the solutions if the dimension $n\le4$.  Typical example in applications are the generalized SKT models (see \cite{SKT}) consisting of more than 2 equations and allowing arbitrary growth conditions in the diffusion and reaction terms (see \refcoro{nis34}). For $n=2$ our \refcoro{nis2} generalizes a result of \cite{NSver} where $A(u,Du)$ was assumed to be independent of $u$.

Next, we will discuss the existence of nontrivial solutions in \refsec{lin}. We now see that \reftheo{extthm} establishes the existence of a strong solution in $\mX$ to \mref{e10}. However, this result provides no interesting information if some 'trivial' or 'semi trivial' solutions, which are solutions to a subsystem of \mref{e10}, are obviously guaranteed by other means. We will be interested in finding other nontrivial solutions to \mref{e10} and the uniform estimates in \refsec{extsec} still play a crucial role here. Although many results in this section, in particular the abstract results in \refsec{ssindex}, can apply to the general \mref{e10} we restrict ourselves to the system
\beqno{b100}\left\{\barr{ll} -\Div(A(u)Du)=\hat{f}(u)& \mbox{in $\Og$},\\
\mbox{Homogenenous Dirichlet or Neumann boundary conditions}& \mbox{on $\partial \Og$}. \earr\right.\eeq 

This problem is the prototype of a general class of nonlinear elliptic systems which arise in numerous applications, where $u$ usually denotes population/chemical density vector of species/agents. Therefore, we will also be interested in finding {\em positive} solutions of this system, i.e. those are in the positive cone $$\mP:=\{u\in\mX\,:\, u=(u_1,\ldots,u_m),\; u_i(x)\ge0 \;\forall x\in\Og\}.$$ Under suitable assumptions on $\hat{f}$, we will show that $T$ can be defined as a map on a bounded set of $\mP$ into $\mP$, i.e. $T$ is a positive map.

If $\hat{f}(0)=0$ then \mref{b100} has the
{\em trivial solution} $u=0$. A solution $u$ is a {\it semi trivial} solution if some components of $u$ are zero. Roughly speaking, we decompose $\mX=\mX_1\oplus\mX_2$, accordingly $\mP=\mP_1\oplus\mP_2$ with $\mP_i$ being the positive cone of $\mX_i$, and write an element of $\mX$ as $(u,v)$ with $u\in\mX_1$, $v\in \mX_2$. Then $w=(u,0)$, with $u>0$, is a semi trivial positive fixed point if $w$ is a fixed point of $T$ in $\mP$ and $u$ is a fixed point of $T|_{\mP_1}$, the restriction of $T$ to $\mP_1$. We then show that the local indices of $T$ at these semi fixed points are solely determined by those of $T|_{\mP_1}$ at $\mP_2$-{\em stable} fixed points.

The existence of nontrivial solutions then follows if the sum of the local fixed point indices at trivial and semi trivial solutions does not at up to the fixed point index of $T$ in $\mP$. Several results on the structure of \mref{e10} will be given to show that this will be the case.

Finally, in \refsec{nontrivial} if Neumann boundary conditions are considered then it could happen that a nontrivial and {\em constant} solution of \mref{b100} exists and solves $\hat{f}(u)=0$. In this case, the conclusion in the previous section does not provide useful information. We are then interested in finding nontrivial {\em nonconstant} solutions to \mref{b100}. The results in this section greatly improve those in \cite{louni,DLT}, which dealt only with systems of two equations, and establish the effect of cross diffusions in 'pattern formation' problems in mathematical biology and chemistry. Besides the fact that our results here can be used for large systems, the analysis provides a systematic way to study pattern formation problems. Further studies and examples will be reported in our forthcoming paper \cite{dlepatt}. We conclude the paper by presenting a simple proof of the fact that nonconstant solutions do not exist if the diffusion is sufficiently large.

\section{A-priori estimates in $W^{1,p}(\Og)$  and H\"older continuity}\eqnoset\label{estsec}

In this section we will establish key estimates for the proof of our main theorem \reftheo{extthm} asserting the existence of strong solutions. Throughout this section, we consider two vector valued functions $U,W$  from $\Og$ into $\RR^m$ and solve
the following system
\beqno{ga1} -\Div(A(W,DU))=\hat{f}(W,DU).\eeq

We will consider the following assumptions on $U,W$ in \mref{ga1}.
\bdes
\item[U.0)] $A,\hat{f}$ satisfy A),F) and SG) with $u=W$ and $\zeta=DU$.
\item[U.1)] $U\in W^{2,2}(\Og)\cap C^1(\Og)$ and $W\in C^{1}(\Og)$. On the boundary $\partial \Og$, $U$ satisfies Neumann or Dirichlet boundary conditions.
\item[U.2)] There is a constant $C$ such that $|DW|\le C|DU|$ a.e. in $\Og$.

\item[U.3)] The following number is finite: \beqno{LLg} \mathbf{\LLg}=\sup_{W\in\RR^m}\frac{|\llg_W(W)|}{\llg(W)}.\eeq
\edes

In the sequel, we will fix a number $q_0>1$ if $n\le4$ and, otherwise,  $q_0>(n-2)/2$ such that
\beqno{pcond} \frac{2q_0-2}{2q_0}=\dg_{q_0}C_*^{-1} \mbox{ for some  $\dg_{q_0}\in(0,1)$.}\eeq Such numbers $q_0, \dg_{q_0}$ always exist if A) and SG) hold. In fact, if $n\le4$ we  choose $q_0>1$ and sufficiently close to 1; if $n>4$, by our assumption SG), we have $\frac{n-4}{n-2}<C_*^{-1}$ and we can choose $q_0>(n-2)/2$ and $q_0$ is sufficiently close to $(n-2)/2$.

The main result of this section shows that if $\|U\|_{BMO(B_R(x_0)\cap\Og)}$ is sufficiently small when $R$ is  uniformly  small then for some $p>n$ $\|DU\|_{L^p(\Og)}$ can be controlled.
\bprop{Dup-uni} Suppose that U.0)-U.3) hold.
Assume that there exists $\mu_0\in(0,1)$, which is sufficiently small, in terms of the constants in A) and F), such that the following holds.
\bdes\item[D)] 
there is a positive $R_{\mu_0}$ such that  \beqno{KeyDumu} \mathbf{\LLg}^2\sup_{x_0\in\bar{\Og}}\|U\|_{BMO(B_{_{\mu_0}}(x_0)\cap\Og)}^2 \le \mu_0 ,\eeq  
\edes

Then   there are $q>n/2$ and a constant $C$ depending on the constants in U.0)-U.3), $q,R_{\mu_0}$, the geometry of $\Og$ and $\|DU\|_{L^2(\Og)}$ such that
\beqno{Viter1zz1} \iidx{\Og}{|DU|^{2q}} \le C.\eeq

In particular, if $U$ is also in $L^1(\Og)$ then $U$ belongs to $C^\ag(\Og)$ for some $\ag>0$ and its norm is bounded by a similar constant $C$ as in \mref{Viter1zz1}.\eprop

The dependence of $C$ in \mref{Viter1zz1} on the geometry of $\Og$ is in the following sense:   Let $\mu_0$ be as in D). We can find balls $B_{R_{\mu_0}}(x_i)$, $x_i\in\bar{\Og}$, such that
\beqno{Ogcover} \bar{\Og} \subset \cup_{i=1}^{N_{\mu_0}}B_{R_{\mu_0}}(x_i),\eeq then $C$ in \mref{Viter1zz1} also depends on the number $N_{\mu_0}$.

The proof of \refprop{Dup-uni} relies on local estimates for the integral of $|DU|$ in finitely many balls $B_{R}(x_i)$ with sufficiently small radius $R$ to be determined by the condition D).
We will establish local estimates for $DU$ in these balls and then add up the results to obtain the global estimate \mref{Viter1zz1}.

\subsection{Local Gagliardo-Nirenberg inequalities involving BMO norms}\label{GNineq}
We first present \reflemm{dleGNlocalz}, one of our main ingredients in the the proof of \refprop{Dup-uni}. This lemma is a simple consequence of  the following local weighted Gagliardo-Nirenberg inequality which is proved in our recent work \cite{dleGNans}. In order to state the assumption for that inequality, we recall some well known notions from Harmonic Analysis.
For $\cg\in(1,\infty)$ we say that a nonnegative locally integrable function $w$ belongs to the class $A_\cg$ or $w$ is an $A_\cg$ weight if the quantity
\beqno{aweight} [w]_{\cg} := \sup_{B_R(y)\subset\Og} \left(\mitx{B_R(y)}{w}\right) \left(\mitx{B_R(y)}{w^{1-\cg'}}\right)^{\cg-1} \; \mbox{ is finite}.\eeq
Here, $\cg'=\cg/(\cg-1)$.  For more details on these classes we refer the readers to \cite{FPW,OP,st}.

We proved in \cite{dleGNans} the following result.
\blemm{dleGNlocal} \cite[Lemma 2.2]{dleGNans} Let $u,U:\Og\to \RR^m$ be vector valued functions with $u\in C^1(\Og)$, $U\in C^2(\Og)$ and $\Fg:\RR^m\to\RR$ be a $C^1$ function such that 
\bdes\item[GN)] $\Fg(u)^\frac{2}{p+2}$ belongs to the $A_{\frac{p}{p+2}+1}$ class. \edes

For any ball $B_t$ in $\Og$
we set \beqno{dleIdeft} I_1(t):=\iidx{B_t}{\Fg^2(u)|DU|^{2p+2}},\,\hat{I}_1(t):=\iidx{B_t}{\Fg^2(u)|Du|^{2p+2}},\eeq
\beqno{dleIdef1t} \bar{I}_1(t):=\iidx{B_t}{|\Fg_u(u)|^2(|DU|^{2p+2}+|Du|^{2p+2})},\eeq and \beqno{dleIdef2t}
I_2(t):=\iidx{B_t}{\Fg^2(u)|DU|^{2p-2}|D^2U|^2}.\eeq

Consider any ball $B_s$ concentric with $B_t$, $0<s<t$, and any nonnegative $C^1$ function $\psi$  such that $\psi=1$ in $B_s$ and $\psi=0$ outside $B_t$. Then, for any $\eg>0$ there are positive constants $C_\eg,C_{\eg,\Fg}$,  depending on $\eg$ and $[\Fg^\frac{2}{p+2}(u)]_{\frac{p}{p+2}+1}$, such that \beqno{dleGNlocalest}\barr{lll}I_1(s)&\le& \eg[I_1(t)+\hat{I}_1(t)]+C_{\eg,\Fg}\|U\|^2_{BMO(B_t)}
\left[\bar{I}_1(t) +  I_2(t)\right]\\ &&+C_\eg\|U\|_{BMO(B_t)}\sup_{x\in B_t}|D\psi(x)|^2\iidx{B_t}{|\Fg|^{2}(u)|DU|^{2p}}.\earr\eeq
\elemm

\brem{GNrem1} By approximation, see \cite{SR}, the lemma also holds for $u\in W^{1,2}(\Og)$ and $U\in W^{2,2}(\Og)$ provided that the quantities $I_1,I_2$ and $\hat{I}_1$ defined in \mref{dleIdeft}-\mref{dleIdef2t} are finite.  \erem

If $\Fg\equiv 1$ then $\bar{I}_1\equiv0$ and we can take $u=U$. The condition GN) is clearly satisfied as $[\Fg]_\cg=1$ for all $\cg>1$ (see \mref{aweight}) and $\Fg_u=0$. It is then clear that we have the following special version of the above lemma.
\blemm{dleGNlocalz}  Let $U:\Og\to \RR^m$ be a vector valued function in $ C^2(\Og)$. For any ball $B_t$ in $\Og$
we set \beqno{dleIdeftz} I_1(t):=\iidx{B_t}{|DU|^{2p+2}},\,
I_2(t):=\iidx{B_t}{|DU|^{2p-2}|D^2U|^2}.\eeq

Consider any ball $B_s$ concentric with $B_t$, $0<s<t$, and any nonnegative $C^1$ function $\psi$  such that $\psi=1$ in $B_s$ and $\psi=0$ outside $B_t$. Then, for any $\eg>0$ there is a positive constant $C_\eg$ such that \beqno{dleGNlocalestz}\barr{lll}I_1(s)&\le& \eg I_1(t)+C_{\eg}\|U\|^2_{BMO(B_t)}
I_2(t)\\ &&+C_\eg\|U\|_{BMO(B_t)}\sup_{x\in B_t}|D\psi(x)|^2\iidx{B_t}{|DU|^{2p}}.\earr\eeq
\elemm

\subsection{The proof of \refprop{Dup-uni}}

In the proof, we will only consider the case when $B_{R}(x_i)\subset\Og$. The boundary case ($x_i\in\partial\Og$) is similar, using the fact that $\partial\Og$ is smooth and a reflection argument to extend the function $U$ or $DU$  outside $\Og$, see \refrem{Nboundaryrem} and \refrem{Dboundaryremz}.

For any  $x_0\in\bar{\Og}$ and $t>0$ we will denote $B_t(x_0)=B_t(x_0)\cap\Og$. For $q\ge1$ we introduce the following quantities. 
\beqno{ABHdef2}\ccB_q(t,x_0)= \iidx{B_t(x_0)}{|DU|^{2q+2}},\eeq 
\beqno{ABHdef2aa} \ccH_q(t,x_0) = \iidx{B_t(x_0)}{|DU|^{2q-2}|D^2U|^2},\eeq \beqno{ABHdef3aa}\ccG_q(t,x_0)= \iidx{B_t(x_0)}{|DU|^{2q}}.\eeq

In the rest of this section, let us fix a point $x_0$ in $\Og$ and drop $x_0$ in the notations \mref{ABHdef2}-\mref{ABHdef3aa}. 

For any $s,t$ such that $0< s<t \le R$ let $\psi$ be a cutoff function for two balls $B_s,B_t$ centered at $x_0$. That is, $\psi$ is nonnegative, $\psi\equiv1$ in $B_s$ and $\psi\equiv0$ outside $B_t$ with $|D\psi|\le1/(t-s)$.

We first have the following local energy estimate.
\blemm{locest}Asume U.0)-U.3). Assume that $q\ge 1$ satisfies the condition \beqno{pcondz} \frac{2q-2}{2q}=\dg_{q}C_*^{-1} \mbox{ for some  $\dg_{q}\in(0,1)$.}\eeq
There is a constant $C_1(q)$ depending also on the constants in A) and F) such that 
\beqno{keyppp}  \ccH_q(s) \le C_1(q)\left[\mathbf{\LLg}^2\ccB_q(t) + \frac{1}{(t-s)^2}\ccG_q(t) \right] \quad 0<s<t\le R.\eeq
\elemm

\bproof 
By the assumption U.1), we can formally differentiate \mref{ga1} with respect to $x$,  more precisely we can use difference quotients (see \refrem{dqrem}), to get the weak form of
\beqno{ga2} -\Div(A_\zeta(W,DU)D^2U+A_W(W,DU)DWDU)=D\hat{f}(W,DU).\eeq

We denote $\bg(W)=\llg^{-1}(W)$. Testing \mref{ga2} with $\fg=\bg(W)|DU|^{2q-2}DU\psi^2$, which is legitimate since $\ccH_q$ is finite, integrating by parts in $x$ and rearranging, we have 
\beqno{energy0}\iidx{\Og}{\myprod{A_\zeta(W,DU)D^2U+A_W(W,DU)DWDU,D\fg}}=\iidx{Q}{\myprod{D\hat{f}(W,DU),\fg}}.\eeq

For simplicity, we will assume in the proof that $\hat{f}\equiv0$. As
$D\fg = I_0+I_1+I_2$ with $$I_0:=\bg(W)D(|DU|^{2q-2}DU)\psi^2,\;I_1:=|DU|^{2q-2}DU\bg_WDW\psi^2,$$ and $\; I_2:=2\bg(W)|DU|^{2q-2}DU\psi D\psi$,
we can rewrite \mref{energy0} as
\beqno{energy0z}\barr{ll}\lefteqn{\iidx{\Og}{\bg(W)\myprod{A_\zeta(W,DU)D^2U,D(|DU|^{2q-2}DU)\psi^2}}}\hspace{.5cm}&\\&=-\iidx{\Og}{[\myprod{A_\zeta(W,DU)D^2U,I_1+I_2}+\myprod{A_W(W,DU)DW,D\fg}]}.\earr\eeq

Let us first 
consider the integral on the left hand side. By U.0) and the uniform ellipticity of $A_\zeta(W,DU)$, we can find a constant $C_*$ such that $|A_\zeta(W,DU)\zeta|\le C_*\llg(W)|\zeta|$. On the other hand, By \mref{pcondz}, $\ag=2q-2$ satisfies $$ \frac{\ag}{2+\ag}=\frac{2q-2}{2q}=\dg_{q}C_*^{-1}=\dg_{q}\frac{\llg(W)}{C_*\llg(W)}.$$
By \cite[Lemma 2.1]{sd}, or \cite[Lemma 6.2]{letrans}, for such $\ag,q$ there is a positive constant $C(q)$  such that 
\beqno{ellest}\myprod{A_\zeta(W,DU)D^2U,D(|DU|^{2q-2}DU)} \ge C(q)\llg(W)|DU|^{2q-2}|D^2U|^2.\eeq

Because $\bg(W)\llg(W)=1$, we then obtain from \mref{energy0z}\beqno{energy1}
\barr{ll}\lefteqn{C_0(q)\iidx{Q}{|DU|^{2q-2}|D^2U|^2\psi^2}}\hspace{.5cm}&\\&\le-\iidx{\Og}{[\myprod{A_\zeta(W,DU)D^2U,I_1+I_2}+\myprod{A_W(W,DU)DW,D\fg}]}.\earr\eeq

The terms $I_1,I_2$ in the integrands on the right hand side of \mref{energy1}  can be easily handled by using the fact that $|D\psi|\le 1/(t-s)$ and the assumption A) which  gives $$|A_\zeta(W,DU)|\le C|\llg(W)| \mbox{ and }|A_W(W,DU)|\le C|\llg_W(W)||DU|.$$ We also note that $|\bg_W(W)|=\llg^{-2}(W)|\llg_W(W)|\le \llg^{-1}(W)\mathbf{\LLg}$ (see \mref{LLg}).

Concerning the first integrand on the right of \mref{energy0z}, using the definition of $I_i$ and Young's inequality, for any $\eg>0$ we can find a constant $C(\eg)$ such that
$$|\myprod{A_\zeta(W,DU)D^2U,I_1}|\le \eg|DU|^{2q-2}|D^2U|^2\psi^2 +C(\eg) \mathbf{\LLg}^2|DW|^2|DU|^{2q}\psi^2,$$
$$|\myprod{A_\zeta(W,DU)D^2U,I_2}|\le \eg|DU|^{2q-2}|D^2U|^2\psi^2 + C(\eg)|DU|^{2q}|D\psi|^2.$$

Similarly, for the second integrand on the right of \mref{energy0z}  we have $$|\myprod{A_W(W,DU)DW,I_0}|\le \eg|DU|^{2q-2}|D^2U|^2\psi^2 + C(\eg)\mathbf{\LLg}^2|DW|^2|DU|^{2q}\psi^2,$$
$$|\myprod{A_W(W,DU)DW,I_1}|\le C\mathbf{\LLg}^2|DW|^2|DU|^{2q}\psi^2,$$
$$|\myprod{A_W(W,DU)DW,I_2}|\le  C\mathbf{\LLg}^2|DW|^2|DU|^{2q}\psi^2+C|DU|^{2q}|D\psi|^2.$$

Choosing $\eg$ sufficiently small, we then obtain from the above inequalities and the assumption $|DW|\le C|DU|$ that
$$\iidx{B_s}{|DU|^{2q-2}|D^2U|^2} \le C_1\mathbf{\LLg}^2\iidx{B_t}{|DU|^{2q+2}}+C_1\frac{1}{(t-s)^2}\iidx{B_t}{|DU|^{2q}}.$$ 
From the notations \mref{ABHdef2}-\mref{ABHdef3aa}, the above estimate gives the lemma. \eproof

\brem{dqrem} For $i=1,\ldots,n$ and $h\ne0$ we denote by $\dg_{i,h}$ the difference quotient operator $\dg_{i,h}u=h^{-1}(u(x+he_i)-u(x))$, with $e_i$ being the unit vector of the $i$-th axis in $\RR^n$. We then apply $\dg_{i,h}$ to the system for $U$ and then test the result with $|\dg_{i,h}U|^{2q-2}\dg_{i,h}U\psi^2$. The proof then continues to give the desired energy estimate by letting $h$ tend to $0$. \erem

\brem{fdurem1} If $\hat{f}\ne0$ then there is an  extra term $|D\hat{f}(W,DU)||DU|^{2q-1}\psi^2$ in \mref{energy1}. This term will give rise to similar terms  in the proof. Indeed, by \mref{FUDU11a} in F) with $u=W$ and $p=DU$, $$|D\hat{f}(W,DU)| \le  C[\llg(W)|D^2U|+ |\llg_W(W)||DW||DU|+ \llg(W)|DU|].$$ As $\bg(W)=\llg^{-1}(W)$ and $|DW|\le C|DU|$, for any $\eg>0$ we can use Young's inequality  and the definition of $\mathbf{\LLg}$ to find a constant $C(\eg)$ such that
$$\barr{lll}|D\hat{f}(W,DU)|\bg(W)|DU|^{2q-1} &\le& C[|D^2U||DU|^{2q-1}+ \mathbf{\LLg}|DU|^{2q+1}+ |DU|^{2q}]\\& \le&
\eg|DU|^{2q-2}|D^2U|^2 +\mathbf{\LLg}^2|DU|^{2q+2}+ C(\eg)|DU|^{2q} .\earr$$

Choosing $\eg>0$ sufficiently small, we then see that the proof can continue to obtain the energy estimate \mref{keyppp}.

Similarly, if we can allow $\hat{f}$ to have nonlinear growth in $DU$ by replacing \mref{FUDU11a} with $$|D\hat{f}(u,p)| \le C[\llg(u)|p|^{\ag-1}|Dp| + |\llg_u(u)||Du||p|^\ag+\llg(u)|p|^\ag] \quad \mbox{for some $\ag\in[1,2)$}.$$ Then $|D\hat{f}(W,DU)|\bg(W)|DU|^{2q-1}$ can be estimated by 
$$C[|D^2U||DU|^{2q+\ag-2}+ \mathbf{\LLg}|DU|^{2q+\ag}+ |DU|^{2q+\ag-1}]$$
Again, by Young's inequality and $q\ge1$ and $\ag<2$, it is not difficult to see that there is some exponent $\cg>0$ depending on $\ag$ such that the above is bounded by
$$\eg|DU|^{2q-2}|D^2U|^2 +(\mathbf{\LLg}^2+\eg)|DU|^{2q+2}+ C(\eg)|DU|^{2q} +C(\eg)(\mathbf{\LLg}^\cg+1),$$ and the proof can continues. 
\erem

\brem{C1rem} Inspecting our proof here and the proof of \cite[Lemma 6.2]{letrans}, we can see that the constant $C(q)$ in \mref{ellest} is decreasing in $q$  and hence $C_1(q)$ is increasing in $q$. Note also that this is the only place where we need \mref{pcondz}.
\erem

\brem{Nboundaryrem} We discuss the case when the centers of $B_\rg, B_R$ are on the boundary $\partial \Og$. We assume that  $U$ satisfies the Neumann boundary condition  on $\partial \Og$. By flattening the boundary we can assume that $B_R\cap\Og$ is the set $$B^+=\{x\,:\, x=(x_1,\ldots,x_n)\mbox{ with } x_n\ge0 \mbox{ and } |x|<R\}.$$

For any point $x=(x_1,\ldots,x_n)$ we denote by $\bar{x}$ its reflection  across the plane $x_n=0$, i.e.,  $\bar{x}=(x_1,\ldots,-x_n)$. Accordingly, we denote by $B^-$ the reflection of $B^+$. For a function $u$ given on $B_+$ we denote its even reflection by  $\bar{u}(x)=u(\bar{x})$ for $x\in B^-$. We then consider the even extension of $\hat{u}$ in $B=B^+\cup B^-$  
$$\hat{u}(x)=\left\{\barr{ll}u(x) &\mbox{if $x\in B^+$},\\\bar{u}(x)&\mbox{if $x\in B^-$}.\earr \right.$$

With these notations, for $x\in B^+$ we observe that  $\Div_x(D_xU)=\Div_{\bar{x}}(D_{\bar{x}}\bar{U})$ and $D_xWD_xU=D_{\bar{x}}\bar{W}D_{\bar{x}}\bar{U}$. Therefore, it is easy to see that $\hat{U}$ satisfies in $B$ a system similar to the one for $U$ in $B^+$. Thus, the proof can apply to $\hat{U}$ to obtain the same energy estimate near the boundary.
\erem

\brem{Dboundaryremz} For the Dirichlet boundary condition we make use of the odd reflection $\bar{u}(x)=-u(\bar{x})$ and then define $\hat{u}$ as in \refrem{Nboundaryrem}. Since $D_{x_i}U=0$ on $\partial\Og$ if $i\ne n$, we can test the system \mref{ga2}, obtained by differentiating the system of $U$ with respect to $x_i$, with $|D_{x_i}U|^{2q-2}D_{x_i}U\psi^2$ and the proof goes as before because no boundary integral terms appear in the calculation. We need only consider the case $i=n$. We observe that $D_{x_n}\hat{U}$ is the even extension of $D_{x_n}U$ in $B$ therefore $\hat{U}$ satisfies a system similar to \mref{ga2}. The proof then continues.
\erem

Next, let us recall the following elementary iteration result in \cite{dleGNans} (which is a consequence of \cite[Lemma 6.1, p.192]{Gius}). 

\blemm{Giusiter1} Let $F,G,g,h$ be bounded nonnegative functions in the interval $[\rg,R]$ with $g,h$ being increasing. Assume that for $\rg \le s<t\le R$ we have \beqno{x1}F(s) \le \eg_0 [F(t)+G(t)]+[(t-s)^{-\ag} g(t)+h(t)],\eeq
\beqno{x2}G(s) \le C[F(t)+(t-s)^{-\ag} g(t)+h(t)]\eeq with $C\ge0$, $\ag,\eg_0>0$.

If $2C\eg_0<1$ then there is constant $c(C,\ag,\eg_0)$ such that \beqno{xx3}F(s)+G(s) \le c(C,\ag,\eg_0)[(t-s)^{-\ag} g(t)+h(t)]\quad \rg \le s<t\le R.\eeq \elemm

We are now ready to give the proof of the main result, \refprop{Dup-uni}, of this section.

\bproof 
For any given $R_0,\eg>0$, multiplying \mref{dleGNlocalestz} by $\mathbf{\LLg}^2$ and using the notations \mref{ABHdef2}-\mref{ABHdef3aa}, we can find a constant $C_\eg$ such that
\beqno{x1z}\mathbf{\LLg}^2\ccB_q(s)\le \eg \mathbf{\LLg}^2\ccB_q(t)+C_{\eg}\mathbf{\LLg}^2\|U\|^2_{BMO(B_t)}
\ccH_q(t)+C_\eg\|U\|_{BMO(B_t)}\frac{\mathbf{\LLg}^2}{(t-s)^2}\ccG_{q}(t)\eeq
for  all $s,t$ such that $0<s<t\le R_0$.

On the other hand, let $q_0>1$ and satisfies \mref{pcond}. We then have  $$ \frac{2q-2}{2q}<C_*^{-1} \quad \forall q\in[1,q_0].$$
Hence \mref{pcondz} of \reflemm{locest} holds for $q\in[1,q_0]$ and we obtain from \mref{keyppp} that \beqno{x2z}\ccH_{q}(s) \le C_1(q) \mathbf{\LLg}^2\ccB_{q}(t) + \frac{C_1(q)}{(t-s)^2}\ccG_{q}(t),\quad 0<s<t\le R_0.\eeq 

We define $$F(t)=\mathbf{\LLg}^2\ccB_q(t),\, G(t)=\ccH_q(t),\, g(t)=\max\{C_1(q_0),C_\eg\|U\|_{BMO(B_{R_0})}\mathbf{\LLg}^2\}\ccG_q(t),$$ 
$$\eg_0=\max\{\eg,C_{\eg}\mathbf{\LLg}^2\|U\|^2_{BMO(B_{R_0})}\}.$$

It is clear that \mref{x1z},\mref{x2z} respectively imply \mref{x1} and \mref{x2} of \reflemm{Giusiter1} with $C=C_1(q_0)$, using the fact that
(see \refrem{C1rem}) $C_1(q)$ is increasing in $q$.

We first choose $\eg$ such that $2C_1(q_0)\eg<1$ and then $R_0>0$ such that \beqno{llg0large}2C_1(q_0)C_{\eg}\mathbf{\LLg}^2\|U\|^2_{BMO(B_{R_0})}<1.\eeq We thus have $2C_1(q)\eg_0<1$ so that \mref{xx3} of \reflemm{Giusiter1} provide a constant $C_2$ depending on $C_1(q_0),\eg_0$ such that
$$\ccH_{q}(s)+\mathbf{\LLg}^2\ccB_q(s) \le \frac{C_2}{(t-s)^2}\ccG_{q}(t),\quad 0<s<t\le R_0.$$

For $t=2s$ the above gives (if $q$ satisfies \mref{pcondz}) \beqno{H1}\ccH_{q}(s)+\mathbf{\LLg}^2\ccB_{q}(s) \le \frac{C_3}{s^2}\iidx{Q_{2s}}{|DU|^{2q}} \quad 0<s\le\frac{R_0}{2}.\eeq

Using this estimate for $\ccB_q(t)$ in \mref{keyppp}, with $s=R_0/4$ and $t=R_0/2$ respectively, we derive
\beqno{H1key}\ccH_{q}(R_1) \le \frac{C_4}{R_1^2}\iidx{B_{R_1}}{|DU|^{2q}},\quad R_1=\frac{R_0}{4}.\eeq

Now, we will argue by induction to obtain a bound for $\ccA_q$ for some $q>n/2$. Suppose that for some $q\ge1$ and $q$ satisfies \mref{pcondz} we can find a constant $C_q$  such that
\beqno{qiter} \iidx{\Og}{|DU|^{2q}} \le C_q,\eeq and that \mref{llg0large} holds then \mref{H1key} implies similar bound for $\ccH_q(R_1)$. We now can cover $\Og$ by $N_{R_1}$ balls $B_{R_1}$, see \mref{Ogcover}, and add up the estimate \mref{H1key} for  $\ccH_q(R_1)$ to obtain a constant $C(\Og,R_1,N_{R_1},q)$ such that
\beqno{DUq}\iidx{\Og}{|DU|^{2q-2}|D^2U|^2} \le C(\Og,R_1,N_{R_1},q,C_q).\eeq

Hence, \mref{qiter} and \mref{DUq} yield another constant $C(\Og,R_1,N_{R_1},q,,C_q)$ such that 
$$ \iidx{\Og}{|DU|^{2q}}+\iidx{\Og}{|DU|^{2q-2}|D^2U|^2} \le C(\Og,R_1,N_{R_1},q,,C_q).$$
Therefore, the $W^{1,2}(\Og)$ norm of $|DU|^q$ is bounded. Let $n_*=n/(n-2)$ (or any number greater than 1 if $n=2$). By Sobolev's inequality, the above implies that there is a constant $C(\Og,R_1,q_*)$ such that  \beqno{qiterzzz}\iidx{\Og}{|DU|^{2qq_*}}\le C(\Og,R_1,N_{R_1},q_*)\quad \mbox{for any $q_*\in(1,n_*]$}.\eeq

We now see that \mref{qiter} holds again with the exponent $q$ being replaced by $qq_*$.

Of course, \mref{qiter} holds for $q=1$ with $C_1=\|DU\|_{L^2(\Og)}^2$. Hence, for suitable choice of an integer $k_0$ and $q_*\in (1,n_*)$ to be determined later we define $L_k= q_*^k$ and repeat
the above argument, with $q=L_k$, $k_0$ times as long as $L_k\le q_0$, $k=0,1,2,\ldots,k_0$. We then obtain from \mref{qiterzzz}
\beqno{qiterzzzz} \iidx{\Og}{|DU|^{2L_{k_0}n_*}} \le C(\|DU\|_{L^2(\Og)},\Og,R_1,N_{R_1},q_*,k_0).\eeq

We now determine $q_*$ and $k_0$. If $n\le4$ we let $k_0=1$ and $q_*=\min\{q_0,n_*\}$.  Otherwise, if $n>4$, it is clear that we can find $q_*\in(1,n_*)$ and an integer $k_0$ such that $L_{k_0}=q_*^{k_0} =q_0$. Since $q_0>1$ if $n\le4$ and $q_0>(n-2)/2$ otherwise, it is clear that $2L_{k_0}n_*>n$ in both cases. 

Therefore, \mref{qiterzzzz} shows that \mref{Viter1zz1} holds for $q=L_{k_0}n_*>n/2$. Since $q_*,k_0$ depend on $q_0$, the constant in \mref{qiterzzzz} essentially depends on the parameters in U.0)-U.3) and the geometry of $\Og$. The proof is complete.
\eproof

\brem{UW21rem} From \refrem{GNrem1}, 
we can see that the conclusion of \refprop{Dup-uni} continues to hold for $U\in W^{2,2}$ as long as  the quantities \mref{ABHdef2}-\mref{ABHdef3aa} are finite for $q\in[1,q_0]$, $q_0$ is fixed in \mref{pcond}.
\erem

\section{Existence of Strong Solutions}\eqnoset\label{extsec}

We now consider the system \mref{e10} in this section. Recall that
\beqno{e1} -\Div(A(u,Du))=\hat{f}(u,Du)\eeq in $\Og$ and $u$ satisfies homogeneous Dirichlet or Neumann boundary conditions on $\partial\Og$.
Throughout this section we will assume that $A,\hat{f}$ satisfy A), F) and SG).

To establish the existence of a strong solution, we embed the systems \mref{e1} in the following family of systems with $\sg\in[0,1]$
$$\left\{\barr{l} -\Div(\hat{A}_\sg(U,DU))=\hat{F}_\sg(U,DU) \mbox{ in $\Og$,}\\
\mbox{$U$ satisfies homogeneous Dirichlet or Neumann BC on $\partial \Og$.}
\earr\right.$$

We will introduce a family of maps $T(\sg,\cdot)$, $\sg\in[0,1]$, acting in some suitable Banach space $\mX$ such that their fixed points are strong solutions to the above system.
We then use Leray-Schauder's fixed point index theory to establish the existence of a fixed point of $T(1,\cdot)$, which is a strong solution to \mref{e1}.

Fixing some  $\ag_0\in(0,1)$, we consider the Banach space
\beqno{Xdef}\mX:=C^{1,\ag_0}(\Og) \;\mbox{(resp. $C^{1,\ag_0}(\Og)\cap C_0(\Og)$)}\eeq if Neumann (resp. Dirichlet) boundary conditions are assumed for \mref{e1}.

For each $w\in\ccX$ and $\sg\in[0,1]$, we  define
$$A_\sg(w)=\int_0^1\partial_2A(\sg w,t\sg Dw)dt.$$ Here and in the sequel, we will also use the notations $\partial_1 g(u,\zeta)$, $\partial_2 g(u,\zeta)$ to denote the partial derivatives of a function $g(u,\zeta)$ with respect to its variables $u,\zeta$.

Assume that there is a family of vector valued functions $\hat{f}_\sg(U,\zeta)$ with $\sg\in[0,1]$, $U\in\RR^n$ and $\zeta\in\RR^{nm}$ such that \bdes
\item[f.0)] $\hat{f}_\sg(U,\zeta)$ is continuous in $\sg$ and $C^1$ in $U,\zeta$. \item[f.1)] $\hat{f}_0(U,\zeta)\equiv0$  and $\hat{f}_1(U,\zeta)=\hat{f}(U,\zeta)$ for all $U,\zeta$.\item[f.2)] $\hat{f}_\sg$ satisfies F) uniformly for $\sg\in[0,1]$. That is, there is a constant $C$ such that for $U\in W^{2,2}(\Og)$ and $W=\sg U$ 
$$|D\hat{f}_\sg(U,DU)| \le  C[\llg(W)|D^2U|+ |\llg_W(W)||DW||DU|+ \llg(W)|DU|],$$ a.e. in $\Og$. \edes

Let $K$ be any constant matrix  satisfying \beqno{Kcondgen}\myprod{Ku,u}\ge k|u|^2 \quad \mbox{for some $k>0$ and all $u\in\RR^m$.}\eeq

For a given $w\in\mX$ and $\sg\in[0,1]$ we consider  the following {\em linear} elliptic system for $u$
\beqno{linsys}\left\{\barr{ll}-\Div(A_\sg(w)Du)+Ku+u=\hat{f}_\sg(w,Dw)+Kw+\sg w& \mbox{in $\Og$,} \\ \mbox{ Homogeneous Dirichlet or Neumann boundary conditions}&\mbox{on $\partial \Og$}.\earr\right.\eeq 

From A) and \mref{Kcondgen} we easily see that the system 
$$\left\{\barr{ll}-\Div(A_\sg(w)Du)+Ku+u=0 & \mbox{in $\Og$} \\ \mbox{ Homogeneous Dirichlet or Neumann boundary conditions}&\mbox{on $\partial \Og$}\earr\right.$$
has $u=0$ as the only solution. From the theory of linear elliptic systems with H\"older continuous coefficient, \mref{linsys} has a unique strong solution $u$.  We then define $T(\sg,w)=u$.

As $A$ satisfies A), $A(\sg U,0)=0$. Hence, for $\sg\in(0,1]$ \beqno{sgA}A_\sg(U)DU=\int_0^1\partial_2A(\sg U,t\sg DU)dtDU= \sg^{-1}A(\sg U,\sg DU).\eeq Meanwhile $A_0(U)=\partial_2A(0,0)$.

We now define \beqno{AFdef1}\hat{A}_\sg(U,\zeta)=\sg^{-1}A(\sg U,\sg \zeta)\; \sg\in(0,1], \; \hat{A}_0(U,\zeta)=\partial_2A(0,0)\zeta. \eeq

The fixed points of $T(\sg,\cdot)$, defined by \mref{linsys}  with $\sg\in[0,1]$, are solutions  the following family of systems
\beqno{famsys}\left\{\barr{l} -\Div(\hat{A}_\sg(U,DU))=\hat{f}_\sg(U,DU)+(\sg-1) U\mbox{ in $\Og$,}\\
\mbox{$U$ satisfies homogeneous Dirichlet or Neumann BC on $\partial \Og$.}
\earr\right.\eeq

\brem{fchoicerem} A typical choice of $\hat{f}_\sg$ in applications is $\hat{f}_\sg(U,\zeta)=\hat{f}(\sg U,\sg \zeta)$. It is not difficult to see that  $\hat{f}_\sg(U,\zeta)$ satisfies f.1)-f.2) if  $\hat{f}$ does.
\erem

\subsection{Existence of Strong Solutions:} The main result of this section is the following result.

\btheo{extthm} We assume that $A,\hat{f}_\sg$ satisfy A), f.0)-f.2) and SG) and that  the following number is finite: \beqno{LLgx} \mathbf{\LLg}=\sup_{W\in\RR^m}\frac{|\llg_W(W)|}{\llg(W)}.\eeq
In addition, we assume that the following conditions hold {\em uniformly} for any solution $U$ to \mref{famsys}.

\bdes
\item[U)] There is a constant $C$ such that \beqno{qiterq1zzz} \|U\|_{W^{1,2}(\Og)} \le C.\eeq
\item[M)] 
for any given $\mu_0>0$ there is a positive $R_{\mu_0}$ for which   \beqno{KeyVMO} \mathbf{\LLg}^2\sup_{x_0\in\bar{\Og}}\|U\|_{BMO(B_{x_0})}^2 \le \mu_0.\eeq  
\edes

Then \mref{e1} has at least one strong solution.
\etheo

\bproof We will use Leray-Schauder's fixed point index theory to establish the existence of a fixed point of $T(1,\cdot)$, which is a strong solution to \mref{e1} and the theorem then follows.  To this end, we  will establish the facts.
\bdes
\item[i)] $T(\sg,\cdot):\mX\to\mX$ is compact for $\sg\in(0,1]$.
\item[ii)] $\mbox{ind}(T(0,\cdot), \mathbf{B}, \mX)=1$ (see the definition of indices below).
\item[iii)] A fixed point $u=T(\sg,u)$ is a solution to \mref{famsys}. For $\sg=1$, such fixed points are solutions to \mref{e1}.
\item[iv)] There is $M>0$, independent of $\sg\in[0,1]$ and $K$, such that any fixed point $u^{(\sg)}\in \mX$ of $T(\sg,\cdot)$ satisfies $\|u^{(\sg)}\|_\mX <M$.
\edes  

Once i)-iv) are established, the theorem follows from the Leray-Schauder index theory. Indeed,  we let $\mathbf{B}$ be the ball centered at $0$ with radius $M$ of $\mX$ and consider the Leray-Schauder indices
\beqno{LSindex}\mbox{ind}(T(\sg,\cdot), \mathbf{B}, \mX) \stackrel{def}{=}deg (Id- T(\sg,\cdot), \mathbf{B}, 0),\eeq where the right hand side denote the Leray-Schauder degree with respect to zero of the vector field $Id- T(\sg,\cdot)$. This number is well defined  because $T(\sg,\cdot)$ is compact (by i)) and $Id- T(\sg,\cdot)$ does not have zero on $\partial\mathbf{B}$ (by iv)).

By the homotopy invariance of the indices, $\mbox{ind}(T(\sg,\cdot), \mathbf{B}, \mX)=\mbox{ind}(T(0,\cdot), \mathbf{B}, \mX)$, which is 1 because of ii). Thus, $T(\sg,\cdot)$ has a fixed point in $\mathbf{B}$ for all $\sg\in[0,1]$. Our theorem then follows from iii).

Using regularity properties of solutions to linear elliptic systems with H\"older continuous coefficients, we see that i) holds. The proof of ii) is standard (see \refrem{T0rem} after the proof). Next, iii) follows from the assumption on $\hat{f}_1$ in f.1). 

Finally, the main point of the proof is iv). We have to establish a uniform estimate for the fixed points of $T(\sg,\cdot)$ in $\mX$. To check iv), let  $u^{(\sg)}\in \mX$ be a fixed point of $T(\sg,\cdot)$, $\sg\in[0,1]$. We need only consider the case $\sg>0$. Clearly,  $u^{(\sg)}$ solves 
$$-\Div(A_\sg(u^{(\sg)})Du^{(\sg)})=\hat{f}_\sg(u^{(\sg)},Du^{(\sg)}) +(\sg-1)u^{(\sg)}$$ so that $U=u^{(\sg)}$ is a strong solution of \mref{famsys}. We need to show that $\|U\|_{\mX}$ is uniformly bounded for $\sg\in[0,1]$.

We now denote $W=\sg U$  and   will show that \refprop{Dup-uni} can be applied to the systems \mref{famsys}.  As we assume \mref{LLgx} and  $W=\sg U$, with $u^{(\sg)}\in\mX$ and $U$ is a strong solution,  the conditions U.1) and  U.2)  are clearly verified. 

We will see that U.0) is verified. Firstly, from \mref{sgA} and the assumption that $A$ satisfies A) and   we  will show that 
$\hat{A}_\sg(U,\zeta)$ satisfies A) too. Indeed,
$$\myprod{\hat{A}_\sg(U, \zeta),\zeta}=\myprod{\sg^{-1}A(\sg U,\sg \zeta),\zeta}
=\myprod{\sg^{-2}A(\sg U,\sg \zeta),\sg\zeta}\ge \llg(\sg U)|\zeta|^2,$$ $$\|\hat{A}_\sg(U,\zeta)\|=\sg^{-1}\|A(\sg U,\sg\zeta)\|\le C_*\llg(\sg U)|\zeta|,$$
$$\|\frac{\partial }{\partial U}\hat{A}_\sg(U,\zeta)\|=\|\partial_1A(\sg U,\zeta)\|\le \llg_{\sg U}(\sg U)|\zeta|.$$
Therefore $\hat{A}_\sg$ satisfies A) with $u=\sg U$.

Secondly, from the assumption f.2) on $\hat{f}_\sg( U, \zeta)$, satisfying F) uniformly for $\sg\in[0,1]$, and the fact that $\llg(W)$ is bounded from below we see that the right hand side of \mref{famsys} satisfies F). Thus, U.0) is satisfied for the system \mref{famsys}.

Finally, it is clear that \mref{KeyVMO} in the assumption M) gives the condition D) of \refprop{Dup-uni}.  The assumption \mref{qiterq1zzz} of U) yields that $\|DU\|_{L^2(\Og)}$ is bounded (see also \refrem{U4rem} after the proof).
More importantly, the uniform bound in \mref{KeyVMO} then gives some positive constants $\mu_0,R(\mu_0)$ such that  \refprop{Dup-uni} applies to  $U=u^{(\sg)}, W=\sg U$ and gives a uniform estimate for $\|u^{(\sg)}\|_{W^{1,2q}(\Og)}$ for some $q>n/2$  and $\sg\in[0,1]$. By Sobolev's imbedding theorems this shows that $u^{(\sg)}$ is H\"older continuous with its norm uniformly bounded with respect to $\sg\in[0,1]$. Since $A$ is $C^1$ in $u$, the results in \cite{Gius} then imply that $Du^{(\sg)}\in C^{\ag}(\Og)$ for any $\ag\in(0,1)$ and its norm is uniformly bounded by a constant independent of $\sg,K$.  We then obtain a uniform estimate for $\|u^{(\sg)}\|_{\mX}$ and iv) is verified.

We then see that \mref{e1} has a solution $u$ in $\mX$. Furthermore, \cite[Chapter 10]{Gius} shows that $u$ is a strong solution. The proof is complete. \eproof

\brem{T0rem} The map $T(0,\cdot)$ is defined by the following {\em linear} elliptic system with constant coeffcients ($A_0:=A_0(w)=\partial_2A(0,0)$ and $\hat{f}_0(w,Dw)\equiv0$)
\beqno{linsys0}-\Div(A_0Du)+Ku+u=Kw\eeq with homogeneous Dirichlet or Neumann boundary conditions. We then consider the following family of systems, with the same boundary conditions, for $\tau\in[0,1]$
\beqno{linsys00}-\Div(A_0Du)+Ku+u=\tau Kw\eeq and define the maps $H(\tau,\cdot)$ on $\mX$ by $H(\tau,w)=u$.  The fixed points $u$ of $H(\tau,\cdot)$ satisfy \mref{linsys00} with $u=w$ so that by testing this with $u$ and using A) and \mref{Kcondgen} we easily see that $u=0$. Similarly, $H(0,\cdot)=0$, a constant map. Thus, by homotopy, $\mbox{ind}(H(1,\cdot), \mathbf{B}, \mX)=\mbox{ind}(H(0,\cdot), \mathbf{B}, \mX)=1$.  Obviously, $T(0,\cdot)=H(1,\cdot)$ so that $\mbox{ind}(T(0,\cdot), \mathbf{B}, \mX)=1$.

\erem

\brem{U4rem} In applications, the assumption on the boundedness of $\|U\|_{W^{1,2}(\Og)}$ in U) can be removed if $\llg(u)$ has a polynomial growth in $|u|$ and $\|U\|_{L^1(\Og)}$ is bounded uniformly. We sketch the proof here. We first observe that $\|U\|_{L^q(\Og)}$ is uniformly bounded. In fact, by \cite[Corollary 2.2]{Gius} and then M), there are constants $C_q, C(q,\mu_0)$ such that for $R\le R_{\mu_0}$ \beqno{giusBMO}\left(\frac{1}{|B_R|}\iidx{B_R}{|U-U_R|^q}\right)^\frac1q\le C_q\|U\|_{BMO(B_R)} \le C(q,\mu_0).\eeq We easily deduct from the above estimate  that there is a constant $C$ depending on $\mu_0,R_{\mu_0},\mathbf{\LLg}$ and $\|U\|_{L^1(\Og)}$ such that $\|U\|_{L_q(B_{R_{\mu_0}})}\le C$.

For $W=\sg U$ we now test the system \mref{ga1} with  and  $\psi=(U-U_{2R})\fg^2$, where $\fg$ is a cut off function for $B_R,B_{2R}$ and satisfies $|D\fg|\le CR^{-1}$. We get
$$ \iidx{B_{2R}}{\llg(W)|DU|^2\fg^2} \le C\iidx{B_{2R}}{(|A(W,DU)||U-U_{2R}||D\fg|\fg+|\hat{f}(W,DU)||\psi|)}.$$

Inspired by the condition f.2), if $\llg$ has a polynomial growth then we can assume that $$|\hat{f}(W,DU)|\le C\llg(W)|DU|+C\llg(W)|U|.$$ Thus, we can use Young's inequality to obtain the following Caccioppoli type estimate
\beqno{Cacci}\iidx{B_R}{\llg(W)|DU|^2} \le  C\iidx{B_{2R}}{\left[R^{-2}\llg(W)+\llg(W)|U|\right]|U-U_{2R}|^2}.\eeq Let $R=R_{\mu_0}/2$. If $\llg(W)$ has a polynomial growth in $W$ and $|W|\le|U|$, we can apply Young's inequality to the right hand side to see that it is bounded in terms of $R_{\mu_0}, \|U\|_{L^q(B_{R_{\mu_0}})}^q$ and the constant $C(q,\mu_0)$ in \mref{giusBMO}. Using a finite covering of $\Og$ and the fact that $\llg(W)$ is bounded from below, we add the above inequalities to obtain a uniform bound for $\|DU\|_{L^2(\Og)}$. Hence, the assumption U) can be removed in this case.
\erem

\brem{extremVMO} We applied \refprop{Dup-uni} to strong solution in the space $\mX$ so that $U,DU$ are bounded and the key quantities $\ccB,\ccH$ are finite. However, the bound provided by the proposition did not involve the supremum norms of $U,DU$ but the local BMO norm of $U$ in M) and the constants in A) and F). 
\erem

In applications, the following corollary of the above theorem will be more applicable.

\bcoro{extcoro} The conclusion of \reftheo{extthm} holds true if U) and M) are replaced by the following condition.

\bdes
\item[M')] There is a constant $C$ such that for any solution $U$ to \mref{famsys}. \beqno{qiterq2zzz} \|U\|_{L^1(\Og)},\;\|DU\|_{L^n(\Og)} \le C.\eeq
\edes
\ecoro

\bproof By H\"older's inequality it is clear that M') implies U). To establish the uniform smallness condition M) we can argue by contradiction. We only sketch the idea of the argument here.  If M) is not true then there are sequences of reals $\{\sg_n\}$ in $[0,1]$ and $\{U_n\}$ of solutions of \mref{famsys} converges weakly to some $U$ in $W^{1,2}(\Og)$ and strongly in $L^2(\Og)$ but $\|U_n\|_{BMO(B_{r_n})}>\eg_0$ for some $\eg_0>0$ and a positive sequence $\{r_n\}$ coverging to 0. We then have $\|U_n\|_{BMO(B_R)}$ converge to $\|U\|_{BMO(B_R)}$ for any given $R>0$. Since $DU_n$ is uniformly bounded in $L^n(\Og)$, it is not difficult to see that $DU\in L^n(\Og)$. Hence, by Poincar\'e's inequality and the continuity of the integral of $|DU|^n$, $\|U\|_{BMO(B_R)}$ can be arbitrarily small. Clearly, if $r_n<R$ then $\|U_n\|_{BMO(B_{r_n})} \le \|U_n\|_{BMO(B_R)}$. Choosing $R$ sufficiently small and letting $n$ tend to infinity, $\|U_n\|_{BMO(B_{r_n})}$ can be arbitrarily small. We obtain a contradiction. Hence, M) is true and the proof is complete.
\eproof

\brem{Xuniest}  Consider a family of systems (not necessarily defined as in the proof) $$\left\{\barr{l} -\Div(\hat{A}_\sg(U,DU))=\hat{F}_\sg(U,DU) \mbox{ in $\Og$,}\\
\mbox{$U$ satisfies homogeneous Dirichlet or Neumann BC on $\partial \Og$,}
\earr\right.$$ which satisfies uniformly the assumptions  A), f.0)-f.2) and SG) and that  the number $\mathbf{\LLg}$ in \mref{LLgx} is bounded. If
any strong solutions $U$ of the family satisfies U) and M) (or M')) uniformly then argument in the proof of \reftheo{extthm} shows that there is a constant $C$ depending only on the parameters in A), f.0)-f.2), SG), U),M) and  $\mathbf{\LLg}$ in \mref{LLgx} such that $\|U\|_\mX\le C$. 
\erem

\subsection{Some Examples:}\label{ssstrong}

We now present some examples in applications where \reftheo{extthm} or \refcoro{extcoro} can apply. The main theme in these examples is to establish the uniform bounds \mref{qiterq2zzz} for the norms $\|\cdot\|_{L^1(\Og)}$ and $\|D(\cdot)\|_{L^n(\Og)}$ of solutions to \mref{famsys}. In fact, under suitable assumptions on the structural of \mref{famsys},  we will show that it is sufficient to control $L^1$ norms of the solutions (see \refrem{U4rem}).

For simplicity we will consider only the following quasilinear system
\beqno{exsys}\left\{\barr{l} -\Div(A(u)Du)=f(u) \mbox{ in $\Og$,}\\
\mbox{$u$ satisfies homogeneous Dirichlet or Neumann BC on $\partial \Og$.}
\earr\right.\eeq

Following \refrem{fchoicerem}, we define $\hat{f}_\sg(u,\zeta)=f(\sg u)$. The corresponding version of \mref{famsys} is \beqno{exsysfam}\left\{\barr{l} -\Div(A(\sg u)Du)=f(\sg u) \mbox{ in $\Og$, $\sg\in[0,1]$,}\\
\mbox{$u$ satisfies homogeneous Dirichlet or Neumann BC on $\partial \Og$.}
\earr\right.\eeq

It is clear that \mref{exsys} is \mref{e1} with $A(u,\zeta),\hat{f}(u,\zeta)$ being $A(u)\zeta, f(u)$.  We will assume that these data satisfy A) and F) and that $\llg(u),f(u)$ have comparable polynomial growths.
\bdes\item[G)] Assume that $\llg(u)\sim (1+|u|)^k$ and $|f(u)|\le C|u|^{l+1}+C$ for some $C,k>0$ and $0\le l\le k$.\edes

We first have the following
\blemm{Ducontrollem} Assume G). There is a constant $C$ such that the following holds true for any solution $u$ to \mref{exsysfam}. \beqno{du2}\iidx{\Og}{(1+|\sg u|^k)|D(\sg u)|^2}\le C\sg^{k+3}
\|u\|_{L^1(\Og)}^{k+2} +C\sg^2\|u\|_{L^1(\Og)}.\eeq

\elemm

\bproof Testing the system \mref{exsysfam} with $\sg^2 u$ and using the ellipticity assumption, we obtain 
$$\sg^2 \iidx{\Og}{\llg(\sg u)|Du|^2}\le C\sg\iidx{\Og}{\myprod{f(\sg u),\sg u}}.$$

From the growth assumptions on $\llg(u), f(u)$ in G) and a simple use of Young's inequality applying to the right hand side of the above inequality, one gets 
\beqno{du2zzz}\iidx{\Og}{(1+|\sg u|^k)|D(\sg u)|^2}\le C\sg \iidx{\Og}{(|\sg u|^{k+2}+|\sg u|)}.\eeq

We now recall the following inequality, which can be proved easily by using a contradiction argument and the fact that $W^{1,2}(\Og)$ is embedded compactly in $L^{2}(\Og)$: For any $w\in W^{1,2}(\Og)$, $\eg>0$ and $\ag\in(0,1]$ there exists a constant $C(\eg,\ag)$ such that
\beqno{ePineq}\iidx{\Og}{|w|^2} \le \eg\iidx{\Og}{|Dw|^2} +C(\eg,\ag)\left(\iidx{\Og}{|w|^\ag}\right)^\frac{2}{\ag}.\eeq

Setting $w=|\sg u|^{\frac{k+2}{2}}$ and noting that $w^2=|\sg u|^{k+2}$ and $|Dw|^2\sim |\sg u|^k|D (\sg u)|^2$. Using the above inequality for $\ag=2/(k+2)$ and sufficiently small $\eg>0$, we deduce from \mref{du2zzz} $$\iidx{\Og}{(1+|\sg u|^k)|D(\sg u)|^2}\le C\sg^{k+3}
\left(\iidx{\Og}{|u|}\right)^{k+2} +C\sg^2\iidx{\Og}{|u|}.$$ This is \mref{du2} and the proof is complete. \eproof

In particular, \mref{du2} implies that 
$$\iidx{\Og}{|Du|^2}\le C\sg^{k+1}
\|u\|_{L^1(\Og)}^{k+2} +C\|u\|_{L^1(\Og)}.$$ Hence, as $k>0$ and $\sg\in[0,1]$, there is a constant $C$ depending only on $\|u\|_{L^1(\Og)}$ such that $\|Du\|_{L^2(\Og)}\le C$ for all solutions to \mref{exsysfam}. The following result immediately follows from this fact and \refcoro{extcoro}.
\bcoro{nis2} Assume G) and that $n=2$. If the solutions of \mref{exsysfam} are uniformly bounded in $L^1(\Og)$ then the system \mref{exsys} has a strong solution. \ecoro

For the case $n=3,4$ we consider some $C^2$ map $P:\RR^m\to\RR^m$ we consider the generalized SKT system \beqno{exsys1}\left\{\barr{l} -\Delta(P(u))=f(u) \mbox{ in $\Og$,}\\
\mbox{$u$ satisfies homogeneous Dirichlet or Neumann BC on $\partial \Og$}.
\earr\right.\eeq This system is a generalized version of the SKT model (see \cite{SKT} where $m=2,n\le2$ and the components of $P(u)$ are assumed to be  quadratics).
The above system is a special case of \mref{exsys} if we set $A(u)=P_u(u)$ and assume A) and F). 

We then have the following
\bcoro{nis34} Assume G) and that $n\le4$. If the solutions of \mref{exsys2} are uniformly bounded in $L^1(\Og)$ then the system \mref{exsys1} has a strong solution. \ecoro

\bproof Since $D(P(\sg u))=\sg A(\sg u)Du$, the system \mref{exsysfam} now reads

\beqno{exsys2}\left\{\barr{l} -\Delta(P(\sg u))=\sg f(\sg u) \mbox{ in $\Og$,}\\
\mbox{$u$ satisfies homogeneous Dirichlet or Neumann BC on $\partial \Og$.}
\earr\right.\eeq If $n=2$ the result was proved in \refcoro{nis2}. We only consider the case $n=4$ as the case $n=3$ is similar and simpler. Again, in this proof, let us denote  $w:=|\sg u|^{\frac{k+2}{2}}$ and $M:=\|u\|_{L^1(\Og)}$. From \mref{du2} we see that $w\in W^{1,2}(\Og)$ and we can find constants $C_i(M)$ such that  \beqno{Pest1}\|w\|_{ W^{1,2}(\Og)}\le \sg C_1(M)\Rightarrow
\|w\|_{ L^{4}(\Og)}\le \sg C_2(M), \eeq using Sobolev's imbedding theorem. From the growth condition on $f$ in G) and Young's inequality,   $|f(\sg u)|\le C(w^2+1)$. Therefore, the above estimates and the equation in \mref{exsys2} imply
\beqno{Pest2}\|f(\sg u)\|_{ L^{2}(\Og)}\le  C_3(M)\Rightarrow
\|\Delta(P(\sg u))\|_{ L^{2}(\Og)}\le \sg C_4(M).\eeq  

On the other hand, since $|D(P(\sg u))|\sim (1+|\sg u|^k)|D(\sg u)|$, we can use H\"older's inequality and \mref{du2} the bound for $\|f(\sg u)\|_{ L^{2}(\Og)}$ to see that
$$\iidx{\Og}{|D(P(\sg u))|}\le \|(1+|\sg u|^k)\|_{L^2(\Og)}\|D(\sg u)\|_{L^2(\Og)} \le \sg C_5(M).  $$
Thus,  $D(P(\sg u))\in L^1(\Og)$.   The last inequality in \mref{Pest2} and Schauder's estimates imply  $\|D^2(P(\sg u))\|_{ L^{2}(\Og)}\le \sg C_6(M)$ for some constant $C_6(M)$. By Sobolev's inequality, $$\|D(P(\sg u))\|_{ L^{4}(\Og)}\le \sg C_7(M).$$ 

Because $A(u)=P_u(u)$ and $D(P(\sg u))=\sg A(\sg u)Du$, we have $Du=\sg^{-1}A^{-1}(\sg u)D(P(\sg u))$. As $A(u)$ is elliptic, its inverse is bounded by some constant $C$. We derive from these facts and the above estimate that  $$\|Du\|_{L^4(\Og)} \le \sg^{-1}C\|D(P(\sg u))\|_{ L^{4}(\Og)}\le  CC_7(M).$$ This gives a uniform estimate for $\|Du\|_{L^4(\Og)}$ and completes the proof. \eproof

We end this section by discussing some special cases where the $L^1$ norm can actually be controlled uniformly so that the above corollaries are applicable. 

Inspired by the SKT model in \cite{SKT} with competitive Lotka-Volterra reactions, we consider the following situation. \bdes\item[SKT)] For some $k>0$ assume that $\llg(u)\sim (1+|u|)^k$ and $f_i(u)=u_i(d_i -g_i(u))$ for some $C^1$ function $g(u)=(g_1(u),\ldots,g_m(u))$ satisfying
\beqno{ggrowth} |g(u)|,\,|u||\partial_{u}g(u)|\le C|u|^k\eeq for some positive constant $C$.\edes

\bcoro{uL1bound} Assume SKT). Suppose that there is a positive constants $C_1$  such that \beqno{gpos1}\sum_i
\myprod{u_ig_i(u),u_i} \ge C_1|u|^{k+2}.\eeq Then there is a strong solution to \mref{exsys} (resp. \mref{exsys1} when $n=2$ (resp. $n\le 4$).

In addition, if Neumann boundary condition is assumed then \mref{gpos1} can be replaced by \beqno{gpos2}\sum_i g_i(u)u_i\ge C_1|u|^{k+1}.\eeq

\ecoro

\bproof We now replace $f(\sg u)$ in \mref{exsysfam} by $f_\sg(u)=(f_{1,\sg}(u),\ldots,f_{m,\sg}(u))$ with $$f_{i,\sg}(u)=\sg^k d_iu_i - u_ig_i(\sg u).$$

Since $|Df_\sg(u)|\le C[\sg^\tau+|g(\sg u)|+ \sg|u||\partial_{\sg u}g(\sg u)|]|Du|$, we see that $f_\sg$ will satisfy f.0)-f.2) if the growth condition \mref{ggrowth} holds. It is easy to see that the argument in the proof of \reflemm{Ducontrollem} continues to hold with this new choice of $f_\sg$ and gives \mref{du2}. Hence, the assertions on existence of strong solutions of the above corollaries continues to hold if we can uniformly control the $L^1(\Og)$ norm of the solutions. This is exactly what we will do in the sequel.

Let us consider the assumption \mref{gpos1} first. We deduce from \mref{gpos1} that
$\myprod{u_ig_i(\sg u),u_i} \ge C_1\sg^k|u|^{k+2}$. Therefore, testing the system \mref{famsys} with $u$, we obtain $$\iidx{\Og}{\llg(\sg u)|Du|^2}\le C_1\sg^{k}\iidx{\Og}{|u|^2}-C_2\sg^{k}\iidx{\Og}{|u|^{k+2}}.$$

Let $w=u$ in \mref{ePineq} and multiply the result with $\sg^k$ to have $$\sg^k\iidx{\Og}{|u|^2} \le \eg\sg^k\iidx{\Og}{|Du|^2} +C(\eg,\ag)\sg^k\left(\iidx{\Og}{|u|}\right)^2.$$

Because $\llg(\sg u)\ge\llg_0>0$, for sufficiently small $\eg$ we deduce from the above two inequalities that there is a constant $C_4$ such that
$$C_2\sg^{k}\iidx{\Og}{|u|^{k+2}}\le C_4\sg^k\left(\iidx{\Og}{|u|}\right)^2.$$

Applying H\"older's inequality to the left hand side integral, we derive
$$C_5\left(\iidx{\Og}{|u|}\right)^{k+2}\le C_4\left(\iidx{\Og}{|u|}\right)^2, \quad C_5>0.$$
Since $k>0$,   the above inequality shows  that $\|u\|_{L^1(\Og)}$ is bounded by a fixed constant.

We now consider the assumption \mref{gpos2} and assume the Neumann boundary condition. Testing the system with $1$, we obtain
$$\iidx{\Og}{g_i(\sg u)u_i}=\sg^k\iidx{\Og}{d_iu_i}.$$

From \mref{gpos2},  $\sum_i g_i(\sg u)u_i\ge C_1\sg^k|u|^{k+1}$. We then derive from the above equation the following
$$C_1\sg^k\iidx{\Og}{|u|^{k+1}}\le C(d_i)\sg^k\iidx{\Og}{|u|}.$$
Again, applying H\"older's inequality to the left hand side integral, we derive
$$C_2\left(\iidx{\Og}{|u|}\right)^{k+1}\le C(d_i)\iidx{\Og}{|u|}$$ for some positive constant $C_2$. Again, as $k>0$, the above gives the desired uniform estimate for $\|u\|_{L^1(\Og)}$. The proof is complete. \eproof

\brem{gposrem} The conditions \mref{gpos1} and \mref{gpos2} on the positive definiteness of $g$ need only be assumed for $u$ such that $|u|\ge M$ for some positive $M$. \erem

\section{On Trivial and Semi Trivial Solutions}\label{lin}\eqnoset
We now see that \reftheo{extthm} establishes the existence of a strong solution in $\mX$ to \mref{e1}. However, the conclusion of this theorem does not provide useful information if some 'trivial' or 'semi trivial' solutions, which are solutions to a subsystem of \mref{e1}, are obviously guaranteed by other means. We will be interested in finding other nontrivial solutions to \mref{e1}. To this end,
we will first investigate these 'trivial' or 'semi trivial' solutions. Several sufficient conditions for nontrivial solutions to exist will be presented in \refsec{ssapp}.

Many results in this section, in particular the abstract results in \refsec{ssindex}, can apply to the general \mref{e1}.
However, for simplicity of our presentation we restrict ourselves to the system
\beqno{b1}\left\{\barr{ll} -\Div(A(u)Du)=\hat{f}(u)& \mbox{in $\Og$},\\
\mbox{Homogenenous Dirichlet or Neumann boundary conditions}& \mbox{on $\partial \Og$}. \earr\right.\eeq 
As in the previous section, we fix some $\ag_0>0$  and let $\mX$ be $C^{1,\ag_0}(\Og,\RR^m)$ (or $C^{1,\ag_0}(\Og)\cap C_0(\Og)$ if Dirichlet boundary conditions are considered). Under appropriate assumptions, \reftheo{extthm} gives the existence of a strong solution in $\mX$ to \mref{b1}. This solution may be trivial. For examples, the 
{\em trivial solution} $u=0$ is a solution to the system if $\hat{f}(0)=0$. 

Let us discuss the existence of {\it semi trivial} solutions. We write $\RR^m=\RR^{m_1}\oplus\RR^{m_2}$ for some $m_1,m_2\ge0$ and denote $\mX_i=C^{1,\ag_0}(\Og,\RR^{m_i})$. By
reordering the equations and variables, we write $\mX=\mX_1\oplus\mX_2$, an element of $\mX$ as $(u,v)$ with $u\in\mX_1$, $v\in \mX_2$,
and  $$A(u,v)=\left[\barr{cc}P^{(u)}(u,v)&P^{(v)}(u,v)\\
Q^{(u)}(u,v)&Q^{(v)}(u,v)\earr\right] \mbox{ and } \hat{f}(u,v)=\left[\barr{c}f^{(u)}(u,v)\\
f^{(v)}(u,v)\earr\right].$$ Here, $P^{(u)}(u,v)$ and $Q^{(v)}(u,v)$ are matrices of sizes $m_1\times m_1$ and $m_2\times m_2$ respectively.

Suppose that \beqno{semicond0}Q^{(u)}(u,0)=0 \mbox{ and } f^{(v)}(u,0)=0 \quad \forall u\in\mX_1,\eeq
then  $(u,0)$, with $u\ne0$, is a semi trivial solution if $u$ solves the subsystem  $$-\Div(P^{(u)}(u,0)Du)=f^{(u)}(u,0).$$

For each $u\in \mX$ and some constant matrix $K$ we consider the following {\em linear} elliptic system for $w$.
\beqno{b2}\left\{\barr{ll} -\Div(A(u)Dw)+Kw=\hat{f}(u)+Ku& \mbox{in $\Og$},\\\mbox{Homogenenous boundary conditions}& \mbox{on $\partial \Og$}.\earr\right.\eeq

For a suitable choice of $K$, see \mref{Kcondgen}, we can always assume that \mref{b2} has a unique weak solution $w\in\mX$. This is equivalent to say that the elliptic system
\beqno{b3}\left\{\barr{ll} -\Div(A(u)Dw)+Kw=0& x\in \Og,\\\mbox{Homogenenous boundary conditions}& \mbox{on $\partial \Og$}\earr\right.\eeq has $w=0$ as the only solution. This is the case if we assume that there is  $k>0$ such that $$\myprod{Ku,u}\ge k|u|^2 \quad \forall u\in\RR^m.$$

We then define $T(u):=w$ with $w$ being the weak solution to \mref{b2}. It is clear that the fixed point solutions of $T(u)=u$ are solutions to \mref{b1}, where $w=u$.

Since $A(u)$ is $C^1$ in $u$, $A(u(x))$ is  H\"older continuous on $\Og$. The regularity theory of linear elliptic systems then shows that $w\in C^{1,\ag}(\Og,\RR^m)$ for all $\ag\in(0,1)$ so that $T$ is compact in $\mX$. Furthermore, $T$ is a differentiable map.

If \mref{b1} satisfies the assumptions of \reftheo{extthm} then there is $M>0$ such that \beqno{Mdef}T(u)=u \Rightarrow \|u\|_\mX<M.\eeq

In applications, we are also interested in finding solutions that are positive.
We then consider  the positive cone in $\mX$
$$\mP := \{u\in\mX\,:\, u=(u_1,\ldots,u_m),\,u_i\ge0\; \forall i\},$$ 
which has nonempty interior
$$\dot{\mP} := \{u\in\mX\,:\, u=(u_1,\ldots,u_m),\,u_i>0\; \forall i\}.$$

Let $M$ be the number provided by \reftheo{extthm} in \mref{Mdef}. We denote by $\mB:=B_\mX(0,M)$ the ball in $\mX$ centered at $0$ with radius $M$. If $T$ maps $\mB\cap\mP$ into $\mP$ then, since $\mP$ is closed in $\mX$ and convex and it is a retract of $\mX$ (see \cite{Dugunji}), we can define the cone index $\mbox{ind}(T,U,\mP)$ for any open subset $U$ of $\mB\cap \mP$ as long as $T$ has no fixed point on $\partial U$, the boundary of $U$ in $\mP$ (\cite[Theorem 11.1]{Amsurv}).

The argument in the proof of \reftheo{extthm} can apply here to give
$$\mbox{ind}(T,\mB\cap \mP,\mP)=1.$$

This yields the existence of a fixed point of $T$, or a solution to \mref{b1}, in $\mP$. From the previous discussion, this solution may be trivial or semi trivial. To establish the existence of a nontrivial positive solution $u$, i.e. $u\in\dot{\mP}$, we will compute the local indices of $T$ at its trivial and semi trivial fixed points. If these indices do not add up to $\mbox{ind}(T,\mB\cap \mP,\mP)=1$ then the existence of nontrivial solutions follows from \cite[Corollary 11.2]{Amsurv}.

\subsection{Some general index results}\label{ssindex}
We then consider the case when \mref{b1} has trivial or semi trivial solutions. That is when $u=0$ or some components of $u$ is zero. We will compute the local indices of the map $T(u)$ at these trivial or semi trivial solutions. The abstract results in this section are in fact independent of \mref{b1} and thus can apply to \mref{e1} and other general situations as well.

We decompse $\mX$ as  $\mX=\mX_1\oplus\mX_2$  and denote by $\mP_i$ and $\dot{\mP}_i$, $i=1,2$, the positive cones and their nonempty interiors in $\mX_i$'s.
We assume (see also \mref{Mdef}) that there is $M>0$ such that the map $T$ is well defined as a map from the ball $\mB$ centered at $0$ with radius $M$  into $\mP$. Accordingly, we denote $\mB_i=\mB\cap\mX_i$.

For $(u,v)\in \mB_1\oplus\mB_2$, we write $$T(u,v)=(F_1(u,v),F_2(u,v)),$$ where $F_i$'s are maps from $\mB$ into $\mX_i$. We also write $\partial_uF_i,\partial_vF_i,...$ for the partial Fr\'echet derivatives of these maps. 

It is clear that for $\fg=(\fg_1,\fg_2)\in \mX_1\oplus\mX_2$ $$T'(u,v)\fg = (\partial_uF_1(u,v)\fg_1+\partial_vF_1(u,v)\fg_2,\partial_uF_2(u,v)\fg_1+\partial_vF_2(u,v)\fg_2).$$

For any fixed $u\in\mB_1$ and $v\in\mB_2$, we will think of  $F_1(\cdot,v)$ and $F_2(u,\cdot)$ as maps from $\mB_1$ into $\mX_1$ and from $\mB_2$ into $\mX_2$ respectively. With a slight abuse of the notation, we still write $\partial_uF_1,\partial_vF_2$ for the Fr\'echet derivatives of these maps.

Taking into account of \mref{semicond0}, we will therefore assume in the sequel that \beqno{F20a} F_2(u,0)=0\quad \forall u\in\mP_1.\eeq This implies
\beqno{F20}F_2(u,tv)=t\int_0^1 \partial_vF_2(u,tsv)v\,ds,\eeq where $\partial_vF_2(u,\cdot)$ is the derivative of $F_2(u,\cdot) : \mB_2\to\mX_2$.

Let $Z_1$ be the set of fixed points of $F_1(\cdot,0)$ in $\mP_1$ and assume that $Z_1\ne\emptyset$. Of course, $u\in Z_1$ iff $F_1(u,0)=u$ and $F_2(u,0)=0$.

For each $u\in \mB_1$ we consider the spectral radius $r_v(u)$ of $\partial_vF_2(u,0)$.
$$r_v(u)=\lim_{k\to\infty}\|\partial_vF_2(u,0)\|_{L(\mX_2)}^{1/k}.$$

We also consider the following subsets of $Z_1$  \beqno{Z12}Z_1^+=\{u\in Z_1\,:\, \mbox{$r_v(u)>1$}\},\,
Z_1^-=\{u\in Z_1\,:\, \mbox{$r_v(u)<1$}\}.\eeq
Roughly speaking, $Z_1^+$ (resp. $Z_1^-$) consists of unstable (resp. stable) fixed points of $T$ in the $\mP_2$-direction. Sometimes we simply say that an element in $Z_1^+$ (resp. $Z_1^-$) is {\em $v$-unstable} (resp. {\em $v$-stable}).

Let us fix an open neighborhood $U$ of $Z_1$ in $\mP_1$. We first need to show that the index $\mbox{ind}(T,U\oplus V)$ is well defined for some appropriate neighborhood of $V$ $0$ in $\mP_2$, i.e. $U\oplus V$ is a neighborhood of $Z_1$ as a subset of $\mP$ and $T$ has no fixed point on its boundary. To this end, we will always assume that \bdes\item[Z)] If $u\in Z_1$ then $\partial_vF_2(u,0)$, the Frechet derivative of $F_2(u,\cdot):\mB_2\to\mX_2$, does not have a positive eigenvector to the eigenvalue 1.\edes

In what follows, if $G$ is a map from an open subset $W$ of $\mP_i$ into $\mP_i$ and there is no ambiguity can arise then we will abbreviate $\mbox{ind}(G,W,\mP_i)$ by $\mbox{ind}(G,W)$. We also say that $G$ is a {\bf strongly} positive endomorphism on $W$ into $\mX_i$ if $G$  maps $W\cap \dot{\mP}_i$ into $\dot{\mP}_i$).

The following main result of this section shows that $\mbox{ind}(T,U\oplus V)$ is determined by the index of the restriction $T|_{\mX_1}$, i.e. $F_1(\cdot,0)$,  at $v$-stable fixed points (in $Z_1^-$).

\btheo{indgen} Assume Z). There is a neighborhood of $V$ of $0$ in $\mP_2$ such that $\mbox{ind}(T,U\oplus V)$ is well defined.  

Suppose also the following. \bdes\item[i)] $T$ is a positive map. That is, $T$ maps $\mB\cap\mP$ into $\mP$.
\item[ii)] $F_2(u,0)=0$ for all $u\in\mB_1$.
\item[iii)] At each semi trivial fixed point $u\in Z_1$, $\partial_v F_2(u,0)$ is a {\bf strongly}  positive map on $\mB_2$ into $\mX_2$. \edes 

Then there exist two disjoint open sets $U^+,U^-$ in $U$ such that $Z_1^+\subset U^+$ and $Z_1^-\subset U^-$ and $$\mbox{ind}(T,U\oplus V)=\mbox{ind}(T,U^-\oplus V)=\mbox{ind}(F_1(\cdot,0),U^-).$$\etheo

\brem{strongpos} For $u\in Z_1$, i) implies that $\partial_v F_2(u,0)$ is a positive endomorphism on $\mX_2$. In fact, for any $u\in Z_1$, $x>0$ and positive small $t$  such that $tx\in V$ we have by our assumptions that $F_2(u,tx)\ge0$ and $F_2(u,0)=0$. Hence, $\partial_v F_2(u,0)x=\lim_{t\to0^+}t^{-1}F_2(u,tx)\ge0$. So that $\partial_v F_2(u,0)$ is positive.  
If certain strong maximum principle for the linear elliptic system defining $\partial_v F_2(u,0)$ is available, see \cite[Theorem 4.2]{Amsurv}, then $\partial_v F_2(u,0)$ is {\em strongly} positive and iii) follows. This assumption can be relaxed if $Z_1$ is a singleton (see \refrem{Irem}). \erem

\brem{posevrem} If $\partial_vF_2(u,0)$ is strongly positive then $r_v(u)$ is the only eigenvalue with positive eigenfunction. Therefore, the assumptions $r_v(u)<1$ and $r_v(u)>1$ are respectively equivalent to the followings \bdes\item[I'.1)] $\partial_vF_2(u,0)$ does not have any positive eigenvector to any eigenvalue $\llg>1$.
\item[I'.2)]  $\partial_vF_2(u,0)$ has a positive eigenvector to some eigenvalue $\llg>1$.
\edes
\erem

The proof of \reftheo{indgen} will be divided into several lemmas which can be of interest in themselves.

Our first lemma shows that there exists a neighborhood $V$ claimed in \reftheo{indgen} such that $\mbox{ind}(T,U\oplus V)$ is well defined.
\blemm{zlemm}
Assume Z). There is $r>0$ such that for $V=B(0,r)\cap \mP_2$, the ball in $\mP_2$ centered at 0 with radius $r>0$, there is no fixed point of $T(u,v)=(u,v)$ with $v>0$ in the closure of $U\oplus V$ in $\mP$. \elemm

\bproof 
By contradiction, suppose that there are sequences $\{r_n\}$ of positives $r_n\to0$ and $\{u_n\}\subset U$, $\{v_n\}\subset \mP_2$ with $\|v_n\|=r_n$ such that, using  \mref{F20}
$$u_n=F_1(u_n,v_n),\, v_n=F_2(u_n,v_n)=\int_0^1\partial_vF_2(u_n,sv_n)v_n\,ds.$$ Setting $w_n=v_n/\|v_n\|$ we have $$w_n=\int_0^1\partial_vF_2(u_n,sr_nw_n)w_n\,ds.$$ 

By compactness, via a subsequence of $\{u_n\}$, and continuity we can let $n\to\infty$ and obtain $u_n\to u$ for some $u\in Z_1$, $v_n\to0$ and $w_n\to w$ in $\mX_2$ such that $u=F(u,0)$ and $\|w\|=1$. Hence,  $w>0$ and satisfies $$w=\int_0^1\partial_vF_2(u,0)w\,ds = \partial_vF_2(u,0)w.$$  Thus, $w$ is a positive eigenvector of $\partial_vF_2(u,0)$   to the eigenvalue 1. This is a contradiction to Z) and completes the proof.
\eproof

In the sequel, we will always denote by $V$ the neighborhood of $0$ in $\mP_2$ as in the above lemma.

Our next lemma on the index of $T$ shows that $T$ can be computed by using its restriction and partial derivatives. 
\blemm{indT} 
We have $$\mbox{ind}(T,U\oplus V)=\mbox{ind}(T_*,U\oplus V), $$ where $T_*(u,v) = (F_1(u,0),\partial_vF_2(u,0)v)$.
\elemm

\bproof Consider the following homotopy
\beqno{homotopy1}H(t,u,v) = \left(F_1(u,tv),\int_0^1 \partial_vF_2(u,tsv)v\,ds\right) \mbox{ for $t\in[0,1]$}.\eeq
We show that this homotopy is well defined on $U\oplus V$. Indeed,  if $H(t,u,v)$ has a fixed point $(u,v)$ on the boundary of $U\oplus V$ for some $t\in[0,1]$ then $$F_1(u,tv)=u,\quad\int_0^1 \partial_vF_2(u,tsv)v\,ds=v,\quad (u,v)\in\partial (U\oplus V).$$ 

Assume first that $t>0$. If $v=0$ then the first equation gives that $F_1(u,0)=u$ so that $u\in Z_1$. But then $(u,0)\notin \partial (U\oplus V)$. Thus, $v>0$ and the second equation (see \mref{F20}) yields $F_2(u,tv)=tv$. This means $(u,tv)$ is a fixed point of $T$ in the closure of $U\oplus V$ with $tv>0$. But there is no such fixed point of $T(u,v)=(u,v)$ in the closure of $U\oplus V$  by \reflemm{zlemm}. Hence, $H(t,u,v)$ cannot have a fixed point $(u,v)$ on the boundary of $U\oplus V$ if $t>0$. We then consider $H(0,u,v)$ whose fixed points $(u,v)\in\partial (U\oplus V)$ must satisfy $u=F_1(u,0)$ so that $u\in Z_1$ and $\partial_vF_2(u,0)v=v$ with $v>0$. But this contradicts Z). 

Thus the homotopy is well defined
and we have that $$\mbox{ind}(T,U\oplus V)=\mbox{ind}(H(1,\cdot),U\oplus V)=\mbox{ind}(H(0,\cdot),U\oplus V).$$  

By \mref{homotopy1}, $H(0,u,v) = (F_1(u,0),\partial_vF_2(u,0)v)=T_*(u,v)$. The proof is complete. \eproof

We now compute $\mbox{ind}(T_*,U\oplus V) $.

\blemm{indT0} Assume that $\partial_v F_2(u,0)$ is a {\bf strongly} positive endomorphism on $\mB_2$ into $\mX_2$ for each $u\in Z_1$.
The following holds

\bdes\item[I.1)] If $r_v(u)<1$ for any $u\in Z_1$  then $\mbox{ind}(T_*,U\oplus V)=\mbox{ind}(F_1(\cdot,0),U)$.
\item[I.2)] If $r_v(u)>1$ for any $u\in Z_1$  then $\mbox{ind}(T_*,U\oplus V)=0$.
\edes
\elemm

\bproof First of all, we see that $\partial_v F_2(u,0)$ is a compact map. In fact, we have  $F_2(u,0)=0$ so that $\partial_v F_2(u,0)x=\lim_{t\to0^+}t^{-1}F_2(u,tx)$.  Since $F_2$ is compact, so is $\partial_v F_2(u,0)$.

To prove I.1), we consider the following hopmotopy $$H(u,v,t) = (F_1(u,0),t\partial_vF_2(u,0)v),\quad t\in[0,1].$$

This homotopy is well defined on $U\oplus V$. Indeed, a fixed point of $(u,v)$ of $H(\cdot,\cdot,t)$ on $\partial(U\oplus V)$ must satisfies $u\in Z_1$ and $tv>0$. But this means $v>0$ is a positive eigenfunction to the eigenvalue $t^{-1}\ge1$. This is a contradiction to Z) and the Krein-Ruthman theorem, see \cite[Theorem 3.2, ii)]{Amsurv} for  strongly positive compact endomorphism on $\mX_2$, $\partial_vF_2(u,0)$ has no positive eigenvector different from $r_v(u)$, which is asumed to be less than 1 in this case. Thus,  
$$\mbox{ind}(T_*,U\oplus V)=\mbox{ind}(H(\cdot,\cdot,0),U\oplus V).$$ But $H(u,v,0)=(F_1(u,0),0)$ so that, by index product theorem, $\mbox{ind}(H(\cdot,\cdot,0),U\oplus V) =\mbox{ind}(F_1(\cdot,0),U)$. Hence,  $\mbox{ind}(T_*,U\oplus V)=\mbox{ind}(F_1(\cdot,0),U)$.

We now consider I.2).  Let $h$ be any element in $\dot{\mP}_2$, the interior of  $\mP_2$. 
We first consider the following homotopy \beqno{homotopy2}H(u,v,t) = (F_1(u,0),t\partial_vF_2(u,0)v+th),\quad t\ge1.\eeq

If $H(\cdot,t)$ has a fixed point $(u,v)$ in $\partial(U\oplus V)$ then $u \in Z_1$ and $v>0$. Thus, there is some $v_*>0$ such that $v_*=t\partial_vF_2(u,0)v_*+th$. This means $t^{-1}v_*-\partial_vF_2(u,0)v_*=h$. Since $t^{-1}\le1<r_v(u)$, this contradicts the following consequence of the Krein-Rutman theorem (see \cite[Theorem 3.2, iv)]{Amsurv} for strongly positive compact operators:
$$\llg x-\partial_vF_2(u,0)x=h \mbox{ has no positive solution if $\llg\le r_v(u)$}.$$

Thus the homotopy is well defined on $U\oplus V$. Because $\partial_vF_2(u,0)v_*\ge0$,  $t\partial_vF_2(u,0)v_*+th$ becomes unbounded as $t\to\infty$, it is clear that $H(u,v,t)$ has no fixed point in $U\oplus V$  for $t$ large. We then have $\mbox{ind}(H(\cdot,\cdot,1),U\oplus V)=0$. 

We now   consider the homotopy $$G(u,v,t)=(F_1(u,0),\partial_vF_2(u,0)v+th) \quad t\in[0,1].$$ 

We will see that this homotopy is well defined on $U\oplus V$ if $\|h\|_{\mX_2}$ is sufficiently small.  First of all, since $\partial_vF_2(u,0)$ is a compact map, the map $f(v)=v-\partial_vF_2(u,0)v$ is a closed map so that $f(\partial V)$ is closed. By Z), if $u\in Z_1$ then $0\not\in f(\partial V)$ so that there is $\eg>0$ such that $B_\eg(0)\cap f(\partial V)=\emptyset$. This means \beqno{vV}\|v-\partial_vF_2(u,0)v\|_{\mX_2}>\eg \quad \forall v\in\partial V.\eeq

We now take $h$ such that $\|h\|_{\mX_2}<\eg/2$. If $G(\cdot,\cdot,t)$ has a fixed point $(u,v)\in\partial(U\oplus V)$ then $u \in Z_1$, $v\in\partial V$ and $v-\partial_vF_2(u,0)v=th$.
This fact and \mref{vV} then yield $$\|v-\partial_vF_2(u,0)v\|_{\mX_2}>\eg>\|th\|_{\mX_2} \quad \forall t\in[0,1].$$ This means  $v-\partial_vF_2(u,0)v \ne th$ for all $u \in Z_1$, $v\in\partial V$. Hence, the homotopy defined by $G$ is well defined on $U\oplus V$.
We then have $$\mbox{ind}(T_*,U\oplus V)=\mbox{ind}(G(\cdot,\cdot,0),U\oplus V)=\mbox{ind}(H(\cdot,\cdot,1),U\oplus V)=0.$$ The proof is complete. \eproof

\brem{Irem} If we drop the assumption that $\partial_v F_2(u,0)$ is strongly positive then the conclusion of \reflemm{indT0} continues to hold if I.1) is replaced by I'.1), which is essentially used in the argument. This is also the case, if we assume I'.2) in place of I.2) and $Z_1$ a singleton, $Z_1=\{u\}$. In fact, let $h$ be such a positive eigenvector of $\partial_vF_2(u,0)$ for some $\llg_u>1$. 
We consider the following homotopy. $$H(u,v,t) = (F_1(u,0),\partial_vF_2(u,0)v+th),\quad t\ge0.$$

We will show that the homotopy is well defined. 
The case $t=0$ is easy. Indeed, if $H(\cdot,\cdot,0)$ has a fixed point $(u,v)$ in $\partial(U\oplus V)$ then $u \in Z_1$ and $\partial_vF_2(u,0)v=v$. But this gives $v=0$, by Z), and $u$ is not in $\partial U$.

We consider the case $t>0$.  
If $H(\cdot,t)$ has a fix point $(u,v)$ in $\partial(U\oplus V)$ then $u \in Z_1$ and $v>0$. Thus, there is some $v_*>0$ such that $v_*=\partial_vF_2(u,0)v_*+th$.  Let $\tau_0$ be the maximal number such that $v_*>\tau_0h$. We then have $\partial_vF_2(u,0)v_*\ge \partial_vF_2(u,0)\tau_0h$ so that (as $\llg>1$) $$v_*=\partial_vF_2(u,0)v_*+th\ge \partial_vF_2(u,0)\tau_0h+th=(\llg\tau_0+t)h>(\tau_0+t)h.$$ Since $t>0$, the above contradicts the maximality of $\tau_0$.
Thus the homotopy is well defined. Again, when $t$ is sufficiently large $H(u,v,t)$ has no solution in $U\oplus V$. Therefore, $\mbox{ind}(H(\cdot,\cdot,1),U\oplus V)=0$. \erem

{\bf Proof of \reftheo{indgen}:} The assumption i) and the regularity results in the previous section show that $T$ is a compact map on $\mX$ so that $\mbox{ind}(T,O,\mP)$ is well defined whenever $T$ has no fixed point on the boundary of an open set $O$ in $\mP$. The assumptions ii) and iii) allow us to make use of the lemmas in this section.

We first prove that $r_v(u)$ is continuous in $u\in Z_1$ (see also \refrem{srcont}). Let $\{u_n\}\subset Z_1$ be a sequence converging to some $u_*\in Z_1$. Accordingly, let $h_n$ be the normalized eigenfunction ($\|h_n\|=1$) to the eigenvalue $\llg_n=r_v(u_n)$. Because $\|\partial_vF_2(u_n,0)\|_{L(\mX_2)}$ is bounded for all $n$,  we see that $\{\llg_n\}$ is bounded from the definition of the spectral radius. Let $\{\llg_{n_k}\}$ be a convergent subsequence of $\{\llg_{n}\}$ converges to some $\llg$. The regularity of elliptic systems yields that the corresponding eigenfunction sequence $\{h_{n_k}\}$ has a convergent subsequence converges to a solution $h>0$ of the eigenvalue problem $\partial_vF_2(u_*,0)h=\llg h$. By uniqueness of the positive eigenfunction (see \cite[Theorem 3.2, ii)]{Amsurv}), $\llg=r_v(u_*)$. We now see that all convergent subsequences of $\{\llg_n\}$ converge to $r_v(u_*)$. Thus, $\limsup \llg_n=\liminf\llg_n$ and $\llg_n=r_v(u_n)\to r_v(u_*)$ as $n\to\infty$. Hence, $r_v(u)$ is continuous in $u\in Z_1$.

Therefore, $Z_1^+,Z_1^-$ are disjoint open sets in $Z_1$. By Z), their union is the compact set $Z_1$ so that they are also closed in $Z_1$ and compact in $\mX_1$. Hence, there are disjoint open sets $U^+,U^-$ in $\mX_1$ such that $Z_1+\subset U^+,Z_1^-\subset U^-$. We then have 
$$\mbox{ind}(T_*,U\oplus V)=\mbox{ind}(T_*,(U^+\cup U^-)\oplus V)=\mbox{ind}(T_*,U^+\oplus V)+\mbox{ind}(T_*,U^-\oplus V).$$

Applying case I.2) of \reflemm{indT0}  for $U=U^+$, we see that  $\mbox{ind}(T_*,U^+\oplus V)$ is zero. It follows that
$$\mbox{ind}(T_*,U\oplus V)=\mbox{ind}(T_*,U^-\oplus V)=\mbox{ind}(F_1(\cdot,0),U^-).$$ By \reflemm{indT}, $\mbox{ind}(T,U\oplus V)=\mbox{ind}(T_*,U\oplus V)$, the theorem then follows. \eproof

\brem{srcont} The strong positiveness of $\partial_vF_2$ is essential in several places of our proof. Under this assumption,  we provided a simple proof of the continuity of $r_v(u)$ on $Z_1$. In general, as $\partial_vF_2$ is always compact, the continuity of $r_v(u)$ follows from \cite[Theorem 2.1]{Gu}, where it was proved that the spectral radius is continuous on the subspace of compact operators. \erem

We end this section by the following well known result which is a special case of \reftheo{indgen}.

\bcoro{indgencoro} Let $X$ be a Banach space with positive cone $P$ and $F$ is a postive compact map on $P$. Suppose that $F(0)=0$ and the directional Frec\'et derivative $F_+'(0)$ exists (i.e. $F_+'(0)x=\lim_{t\to0^+}t^{-1}F(tx)$). Assume also that $F_+'(0)$ does not have any positive eigenvector to the eigenvalue $1$ and that the following holds.
\beqno{Stable}\mbox{$F_+'(0)$ does not have any positive eigenvector to any eigenvalue $\llg>1$}.\eeq Then we can find a neighborhood $V$ of $0$ in $P$ such that $0$ is the only fixed point of $F$ in $V$ and
\beqno{F1ind00} \mbox{ind}(F,V)=\left\{\barr{ll}1 &\mbox{if \mref{Stable} holds},\\ 0 &\mbox{otherwise}. \earr
\right.\eeq \ecoro

To see this,  we let $\mX=\{0\}\oplus X$, i.e. $\mX_1=\{0\}$ and $\mX_2=X$, and $T(\cdot)=(0,F(\cdot))$. Obviously, \reftheo{indgen}, with $F_1$ is the constant map and $F_2=F$, $Z_1$ being the singleton $\{0\}$ (see \refrem{Irem}) and $U=\{0\}$, provides a neighborhood $U^-$ of $Z_1^-$ in $\mX_1$ such that 
$\mbox{ind}(T,\{0\}\oplus V)=\mbox{ind}(0,U^-)$. Clearly, as $F_1$ is a constant map, if \mref{Stable} holds then $U^-=\{0\}$ and $\mbox{ind}(F_1(\cdot,0),U^-)=1$; otherwise $U^-=\emptyset$ and $\mbox{ind}(F_1(\cdot,0),U^-)=0$.
By the product theorem of indices, $\mbox{ind}(T,\{0\}\oplus V)=\mbox{ind}(F,V)$,  \mref{F1ind00} then follows.

\subsection{Applications}\label{ssapp} 

In this section, we will show that the abstract results on the local indices of $T$  at trivial and semi trivial solutions in \reftheo{indgen} can apply  to the map $T$ defined by \mref{b1} satisfying a suitable set of assumptions.

Going back to the definition of $T$,  for each $(u,v)\in \mX$ and some suitable constant matrix $K$ we consider the following {\em linear} elliptic system for $w=T(u,v)$.
\beqno{b2z}\left\{\barr{ll} -\Div(A(u,v)Dw)+Kw=\hat{f}(u,v)+K(u,v)& \mbox{in $\Og$},\\\mbox{Homogenenous boundary conditions for $w$}& \mbox{on $\partial \Og$ }.\earr\right.\eeq

\brem{Kchoice}We observe that the choice of the matrix $K$ is not important here as long as the map $T$ is well defined (as a positive map). In fact, let $K_1,K_2$ be two different matrices and $T_1,T_2$ be the corresponding maps defined by \mref{b2z}. It is clear that these maps have the same set of fixed points consisting of solutions to \mref{b1}. Hence,  via a simple homotopy $tT_1+(1-t)T_2$ for $t\in[0,1]$, the indices $\mbox{ind}(T_i,U)$ are equal whenever one of their indices is defined (i.e. \mref{b1} does not have any solution on $\partial U$).\erem

{\bf Trivial solution:} It is clear that $0$ is a solution if $\hat{f}(0)=0$. In this case, we can apply \refcoro{indgencoro} with $F=T$. The eigenvalue problem of $T'(0)h=\llg h$ now is
\beqno{trivlin} -\Div(A(0)Dh)+ Kh=\llg^{-1}(\hat{f}_u(0)+K)h.\eeq

We then have the following result from \refcoro{indgencoro}.
\blemm{trivlemm} There is a neighborhood $V_0$ of $0$ in $\mP$ such that if \mref{trivlin} has a positive solution $h$ to some eigenvalue $\llg>1$ then $\mbox{ind}(T,V_0)=0$. Otherwise, 
$\mbox{ind}(T,V_0)=1$.\elemm

{\bf Semitrivial solution:} By reordering the equations and variables, we will write an element of $\mX$ as $(u,v)$ 
and  $$A(u,v)=\left[\barr{cc}P^{(u)}(u,v)&P^{(v)}(u,v)\\
Q^{(u)}(u,v)&Q^{(v)}(u,v)\earr\right] \mbox{ and } \hat{f}(u,v)=\left[\barr{c}f^{(u)}(u,v)\\
f^{(v)}(u,v)\earr\right].$$

The existence of semitrivial solutions $(u,0)$ usually comes from the assumption that \beqno{semicond}Q^{(u)}(u,0)=0 \mbox{ and } f^{(v)}(u,0)=0 \quad \forall u\in\mX_1.\eeq

If  \mref{semicond} holds then it is clear that $(u,0)$ is a solution of \mref{b1} if and only if $u$ solves the following subsystem  \beqno{subsysz}-\Div(P^{(u)}(u,0)Du)=f^{(u)}(u,0).\eeq

Let us then assume that the set $Z_1$ of positive solutions to \mref{subsysz} is nonempty.

To compute the local index of $T$ at a semi trivial solution we consider the following matrix \beqno{Kdef}K=\left[\barr{cc}K_1^{(u)}&K_1^{(v)}\\
0&K_2^{(v)}\earr\right],\eeq where the matrices $K_1^{(u)},K_1^{(v)}$ and
$K_2^{(v)}$ are of sizes $m_1\times m_1$, $m_1\times m_2$ and $m_2\times m_2$ respectively. The system in \mref{b2z} for $w=(w_1,w_2)$ now reads
\beqno{b2z1}\barr{ll}\lefteqn{-\Div(P^{(u)}(u,v)Dw_1+P^{(v)}(u,v)Dw_2)+K_1^{(u)}w_1+K_1^{(v)}w_2=}\hspace{7cm}&\\&f^{(u)}(u,v)+K_1^{(u)}u+K_1^{(v)}v,\earr\eeq and
\beqno{b2z2} -\Div(Q^{(u)}(u,v)Dw_1+Q^{(v)}(u,v)Dw_2)+K_2^{(v)}w_2=f^{(v)}(u,v)+K_2^{(v)}v.\eeq

We will consider the following assumptions on the above subsystems.
\bdes\item[K.0)] Assume that there are $k_1,k_2>0$ such that
\beqno{Kcond}  \myprod{K_1^{(u)}x_1,x_1}\ge k_1|x_1|^2, \,\myprod{K_2^{(v)}x_2,x_2}\ge k_2|x_2|^2 \quad \forall x_i\in\RR^{m_i}, i=1,2.\eeq

\item[K.1)] For all $(u,v)\in\mB\cap\mP$
$$\hat{f}(u,v)+K(u,v)\ge0,$$
and the following maximum priciple holds: if $(u,v)\in \mB\cap\mP$ and $w$ solves $$\left\{\barr{ll} -\Div(A(u,v)Dw)+Kw\ge0& \mbox{in $\Og$},\\\mbox{Homogenenous boundary conditions}& \mbox{on $\partial \Og$}\earr\right.$$ then $w\ge0$.
\item[K.2)] For any $u\in Z_1$ and $\fg_2\in\dot{\mP}_2$ a strong maximum principle holds for the system \beqno{hhh1}-\Div(Q^{(v)}(u,0)D\mathbf{U}_2 +Q^{(u)}_{v}(u,0)Du\fg_2)+K_2^{(v)}\mathbf{U}_2=f_{v}^{(v)}(u,0)\fg_2+K_2^{(v)}\fg_2.\eeq That is, if $f_{v}^{(v)}(u,0)\fg_2+K_2^{(v)}\fg_2\in\dot{\mP}_2$ then $\mathbf{U}_2\in\dot{\mP}_2$.
\edes

Concerning the term $Q^{(u)}_{v}(u,0)Du\fg_2$ in K.2) we have used the following notation: if $B(u,v) =(b_{ij}(u,v))$, with $i=1,\ldots,m_2$ and $j=1,\ldots,m_1$,  and $\fg_2=(\fg^{(1)},\ldots,\fg^{(m_2)})$ then \beqno{Bnotation}B_{v}(u,v)Du\fg_2 =\left(\partial_{v^{(k)}}b_{ij}(u,v)\fg^{(k)}Du_j\right)_i=\left(\partial_{v^{(k)}}b_{ij}(u,v)Du_j\right)_{k,i}\fg_2.\eeq 

We also assume that
\bdes \item[K.3)] For any $u\in Z_1$ the linear system $$-\Div(Q^{(v)}(u,0)Dh_2 +Q^{(u)}_{v}(u,0)Duh_2)=f_{v}^{(v)}(u,0)h_2$$ has no positive solution $h_2$ in $\mP_2$.
\edes

For $u\in Z_1$ we will also consider the following eigenvalue problem for an eigenfunction $h$
\beqno{h1z}-\Div(\llg Q^{(v)}(u,0)Dh +Q^{(u)}_{v}(u,0)Duh)+\llg K_2^{(v)}h=f_{v}^{(v)}(u,0)h+K_2^{(v)}h,\eeq and denote
$$Z_1^-:=\{u\in Z_1\,:\, \mbox{\mref{h1z} has a positive solution $h$ to an eigenvalue $\llg<1$}\}.$$

The main theorem of this subsection is the following.

\btheo{indexthm1} Assume K.0)-K.4) with $k_1$ in K.0) being sufficiently large. Then the map $T$ described in \mref{b2z} is well defined on $\mB\cap\mP$ and maps $\mB\cap\mX$ into $\mP$. There are neighborhoods $U,U^-$ respectively of $Z_1,Z_1^-$ in $\mP_1$ and a neighborhood $V$ on $0$ in $\mP_2$ such that $$\mbox{ind}(T,U\oplus V)=\mbox{ind}(F_1(\cdot,0),U^-).$$ Here, $F_1(\cdot,0)$ maps $\mB\cap\mP_1$ into $\mP_1$ and $w_1=F_1(u,0)$, $u\in\mB\cap\mP_1$, is the unique solution to \beqno{ueqn111}-\Div(P^{(u)}(u,0)Dw_1)+K_1^{(u)}w_1=f^{(u)}(u,0) +K_1^{(u)}u.\eeq
\etheo

The above theorem is just a consequence of \reftheo{indgen} applying to the system \mref{b2z}. We need only to verify the assumption of the theorem. For this purpose and later use in the section we will divide its proof into lemmas which also contain additional and useful facts.

We first have the following lemma which shows  that the assumption \mref{F20a}, that $F_2(u,0)=0$ for all $u\in \mX_1$, in the previous section is satisfied.

\blemm{Klemm} Let $T$ be defined by \mref{b2z}. If K.0) holds for some sufficiently large $k_1$ then $T$ is well defined by \mref{b2z} for any given matrix $K_1^{(v)}$. The components $F_1,F_2$ of $T$ satisfy

\bdes\item[i)] $F_2(u,0)=0$ for all $u\in \mX_1$.
\item[ii)] $w_1=F_1(u,0)$ solves \mref{ueqn111}. In addition, $\mbox{ind}(F_1(\cdot,0),\mB\cap\mP_1)=1$.
\edes

\elemm

\bproof  We write $w=T(u,v)=(F_1(u,v),F_2(u,v))$ in \mref{b2z} by $(w_1,w_2)$,  with $w_i\in\mX_i$. Because $$\myprod{Kw,w}=\myprod{K_1^{(u)}w_1,w_1}+ \myprod{K_1^{(v)}w_2,w_1}+\myprod{K_2^{(v)}w_2,w_2},$$  a simple use of Young's inequality and \mref{Kcond} show that if $k_1$ sufficiently large  then $\myprod{Kx,x}\ge |x|^2$ for any given $K_1^{(v)}$,.
Hence, $T$ is well defined by \mref{b2z}.

At $(u,0)$, since  $\hat{f}^{(v)}(u,0)=0$ and $Q^{(u)}(u,0)=0$, the subsystem \mref{b2z2} defining $w_2=F_2(u,0)$ now is
$$-\Div(Q^{(v)}(u,0)Dw_2) +K_2^{(v)}w_2=0.$$  This system has $w_2=0$ as the only solution because of the assumption \mref{Kcond}  on $K_2^{(v)}$ and the ellipticity of $Q^{(v)}(u,0)$. This gives i).

Next, as $w_2$ and $Dw_2$ are zero, \mref{b2z1} gives that $w_1=F_1(u,0)$ solves $$-\Div(P^{(u)}(u,0)Dw_1)+K_1^{(u)}w_1=f^{(u)}(u,0) +K_1^{(u)}u.$$ Again, for a given $u\in\mX_1$ this subsystem has a unique solution $w_1$ if $\myprod{K_1^{(u)}x,x}\ge k_1|x|^2$ for some $k_1>0$. Moreover, the fixed point of $u=F_1(u,0)$ solves $$-\Div(P^{(u)}(u,0)Du)=f^{(u)}(u,0) .$$ This system satisfies the same set of structural conditions for the full system \mref{b1} so that \reftheo{extthm} can apply here to give ii).
\eproof

From the proof of \reftheo{indgen} we need study the Frechet (directional) derivative of $T$ defined by \mref{b2z}. For this purpose and later use,  we consider a more general linear system defining $w=T(u)$
\beqno{b2a}\left\{\barr{ll} -\Div(A(u)Dw+B(u,Du)w)+C(u,Du)w=\hat{f}(u,Du)& x\in \Og,\\\mbox{$w=0$ on $\partial \Og$}, &\earr\right.\eeq for some matrix valued functions $A,B,C,\hat{f}$.

We then recall the following elementary result on the  linearization of the above system at $u$.

\blemm{linlemm} Let $u,\fg$ be in $\mX$. If $w=T(u)$ is defined by \mref{b2a} then  $W=T'(u)\fg$ solves the following system
$$
-\Div(A(u)DW+B(u,Du)W+\ccB(u,w,\fg)) +\ccC(u,W,w,\fg)=\ccF(u,\fg),$$ where
$$\ccB(u,w,\fg)=A_{u}(u)\fg Dw+B_{u}(u,Du)\fg w+B_{\zeta}(u,Du)D\fg w,$$
$$\ccC(u,W,w,\fg)=C(u,Du)W +C_{u}(u,Du)\fg  w+C_{\zeta}(u,Du)D\fg w,$$
$$\ccF(u,\fg)=\hat{f}_u(u,Du)\fg+\hat{f}_\zeta(u,Du)D\fg.$$

\elemm

The proof of this lemma is standard. Because $A,\hat{f}$ are $C^1$ in $u$, it is easy to see that $T$ is differentiable. In fact, for any $u,\fg\in \mX$  we can compute $T'(u)\fg=\lim_{h\to0}\dg_{h,\fg}T(u)$, where $\dg_{h,\fg}$ is  the difference quotient operator  $$\dg_{h,\fg}T(u)= h^{-1}(T(u+h\fg)-T(u)).$$  

Subtracting \mref{b2} at $u$ being $u+h\fg$ and $u$ and dividing the result by $h$, we get
\beqno{difflin}\barr{ll}\lefteqn{-\Div(\dg_{h,\fg}[A(u)DT(u)+B(u,Du)T(u)])+}\hspace{3cm}&\\&
\dg_{h,\fg}[C(u,Du)T(u)]=\dg_{h,\fg}\hat{f}(u,Du).\earr\eeq

It is elementary to see that if $g$ is a $C^1$ function in $u,\zeta=Du,w=T(u),Dw$ then $$\lim_{h\to0}\dg_{h,\fg}g(u,Du,w,Dw) = g_{u}\fg+g_{\zeta}D\fg+g_{w}T'(u)\fg+g_{Dw}D(T'(u)\fg).$$

Using the above in \mref{difflin} and rearranging the terms, we obtain the lemma.

Applying \reflemm{linlemm} to the system \mref{b2z2}, we have the following lemma concerning the map $\partial_vF_2(u,0)$ at $u\in Z_1$. 
\blemm{Tlinlem} Let  $u\in Z_1$. An eigenvector function $h$ of $\partial_vF_2(u,0)h=\llg h$ satisfies the system
\beqno{h1}-\Div(\llg Q^{(v)}(u,0)Dh +Q^{(u)}_{v}(u,0)Duh)+\llg K_2^{(v)}h=f_{v}^{(v)}(u,0)h+K_2^{(v)}h.\eeq

In addition, if K.2) holds then $\partial_vF_2(u,0)$ is strongly positive.
\elemm

\bproof Let $\fg=(0,\fg_2)$ and $u\in Z_1$. We have $$w:=T(u,0),\;W:=T'(u,0)\fg = (\partial_vF_1(u,0)\fg_2,\partial_vF_2(u,0)\fg_2)$$ satisfy, by \reflemm{linlemm} with $B(u,Du)=0$ and $C(u,Du)=K$ 
$$-\Div(A(u,0)DW +\partial_{u,v}A(u,0)\fg Dw )+KW=\hat{f}_{u,v}(u,0)\fg +K\fg.$$

At $(u,0)$, we have that $v,Dv$ are zero and $Dw=D(T(u,0))=(Du,0)$ so that
$$A(u,0)=\left[\barr{cc}P^{(u)}(u,0)&P^{(v)}(u,0)\\
0&Q^{(v)}(u,0)\earr\right],$$
$$\partial_{u,v}A(u,0)\fg Dw=\left[\barr{c}P^{(u)}_v(u,0)Du\fg_2\\
Q^{(u)}_v(u,0)Du\fg_2\earr\right].$$

Thus, $\mathbf{U}_1:=\partial_vF_1(u,0)\fg_2$ and $\mathbf{U}_2:=\partial_vF_2(u,0)\fg_2$, the components of $T'(u,0)\fg$, satisfy
$$\barr{ll}\lefteqn{-\Div(P^{(u)}(u,0)D\mathbf{U}_1+P^{(v)}(u,0)D\mathbf{U}_2 +P^{(u)}_{v}(u,0)Du\fg_2)}\hspace{2cm}&\\&
+K_1^{(u)}\mathbf{U}_1+K_1^{(v)}\mathbf{U}_2=f_{v}^{(u)}(u,0)\fg_2+K_1^{(v)}\fg_2,\earr$$ and
\beqno{h1zz}-\Div(Q^{(v)}(u,0)D\mathbf{U}_2 +Q^{(u)}_{v}(u,0)Du\fg_2)+K_2^{(v)}\mathbf{U}_2=f_{v}^{(v)}(u,0)\fg_2+K_2^{(v)}\fg_2.\eeq

We consider the eigenvalue problem $\partial_vF_2(u,0)h=\llg h$. Set $\fg_2=h$ then $\mathbf{U}_2=\llg h$ and it is clear from \mref{h1zz} that $h$ is the solution to \mref{h1}. 

Finally, the system \mref{h1zz} defining $\mathbf{U}_2:=\partial_vF_2(u,0)\fg_2)$ is exactly \mref{hhh1} in K.2). Thus, the strong maximum principle for \mref{hhh1} yields that $\partial_vF_2(u,0)$ is strongly positive.
\eproof

{\bf Proof of \reftheo{indexthm1}:} \reflemm{Klemm} shows that $T$ is well defined and maps $\mB\cap\mP$ into $\mP$ if K.1) is assumed. \reflemm{Tlinlem} and K.2) then gives the strong positivity of $\partial_v F_2(u,0)$ for any $u\in Z_1$. In addition, the equation in the condition K.3) is \mref{h1} of \reflemm{Tlinlem} when $\llg=1$ so that K.3)  means that the condition Z) of the previous section holds here. Thus, our theorem is just a consequence of \reftheo{indgen}. \eproof

We now turn to semi trivial fixed points of $T$ in $\mX_2$. These fixed points are determined by the following system,  setting $w_1=u=0$, $w_2=v$ in \mref{b2z1} and \mref{b2z2}
\beqno{b2z11}\left\{\barr{l}-\Div(P^{(v)}(0,v)Dv)=f^{(u)}(0,v),\\ -\Div(Q^{(v)}(0,v)Dv)=f^{(v)}(0,v).\earr\right.\eeq

We will assume that this system has no positive solution $v$. In fact, if $P^{(v)}(0,v)\ne0$ the above system is {\em overdetermined} so that the existence of a nonzero solution $v$ of the second subsystem satisfying the first subsystem is very unlikely. In fact, assuming $f^{(v)}(0,0)=0$, it could happen that the second subsystem already has $v=0$ as the only solution.

At $(u,v)=(0,0)$, the eigenvalue problem $\partial_vF_2(0,0)h_2=\llg h_2$ for $h_2\in\mX_2$ is
\beqno{h10}-\Div(Q^{(v)}(0,0)Dh_2)+ K_2^{(v)}h_2=\llg^{-1}(f_{v}^{(v)}(0,0)+K_2^{(v)})h_2.\eeq

From ii) of \reflemm{Klemm}, the eigenvalue problem $\partial_uF_1(0,0)h_1=\llg h_1$ for $h_1\in \mX_1$ is
\beqno{h100}-\Div(P^{(u)}(0,0)Dh_1)+ K_1^{(u)}h_1=\llg^{-1}(f_{u}^{(u)}(0,0)+K_1^{(u)})h_1.\eeq Again, we will say that $0\in\mX_1$ is $u$-stable if the above has no positive eigenvector $h_1$ to any eigenvalue $\llg>1$. Otherwise, we say that $0$ is $u$-unstable.

Our first application of \reftheo{indgen} is to give sufficient conditions such that semi trivial and nontrivial solutions exist.

\btheo{indgenexist} Suppose that \mref{b2z11} has no positive solution. Assume K.0)-K.4) and that the system \beqno{k33} -\Div(P^{(u)}(0,0)Dh_1)=f_{u}^{(u)}(0,0)h_1\; \mbox{ has no solution $h_1\in\dot{\mP}_1$.}
\eeq

If either one of the followings holds
\bdes\item[i.1)] $0$ is  $u$-stable and $Z_1^+=\{0\}$; 
\item[i.2)] $0$ is  $u$-unstable and $Z_1^-=\{0\}$; \edes then there is a nontrivial positive solution to \mref{b1}.
\etheo

\bproof We denote $Z_p=\{u\in Z_1\,:\, u>0\}$. Thus $Z_p$
is the set of semi trivial solutions and $Z_1=\{0\}\cup Z_p$. Accordingly, we denote by $Z_p^+$ (resp. $Z_p^-$) the $v$-unstable (resp. $v$-stable) subset of $Z_p$.
The assumption \mref{k33} means $\partial_u F_1(0,0)$ does not have positive eigenfunction to the eigenvalue 1 in $\mP_1$.  Applying \refcoro{indgencoro} with $X=\mX_1$ and $F(\cdot)=F_1(\cdot,0)$, we can find a neighborhood $U_0$ in $\mX_1$ of $0$ such that $0$ is the only fixed point of $F_1(\cdot,0)$ in $U_0$ and \mref{F1ind00} gives \beqno{F1ind0}\mbox{ind}(F_1(\cdot,0),U_0)=\left\{\barr{ll}1 &\mbox{if $0$ is $u$-stable},\\ 0 &\mbox{if $0$ is $u$-unstable}. \earr
\right.\eeq

Since $Z_1$ is compact and $Z_1=\{0\}\cup Z_p$, the above argument shows that $Z_p$ is compact.   From the proof of \reftheo{indgen}, there are disjoint open neighborhoods $U_p^-$ and $U_p^+$ in $\mX_1$ of $Z_p^-$ and $Z_p^+$ respectively.  Of course, we can assume that $U_0$, $U_p^-$ and $U_p^+$ are disjoint so that for $U=U_0\cup U_p^+\cup U_p^-$ \beqno{indf1}\mbox{ind}(F_1(\cdot,0),U)=\mbox{ind}(F_1(\cdot,0),U_0)+\mbox{ind}(F_1(\cdot,0),U_p^+)+\mbox{ind}(F_1(\cdot,0),U_p^-).\eeq

By \reflemm{Klemm}, $\mbox{ind}(F_1(\cdot,0),\mB\cap\mP_1)=1$. This implies $\mbox{ind}(F_1(\cdot,0),U)=1$. If i.1) holds then $Z_p^+=\emptyset$ and we can take $U_p^+=\emptyset$ and $U^-=U_p^-$ in \reftheo{indgen}. From  \mref{F1ind0}, $\mbox{ind}(F_1(\cdot,0),U_0)=1$ so that \mref{indf1} implies $\mbox{ind}(F_1(\cdot,0),U_p^-)=0$. This yields   $\mbox{ind}(F_1(\cdot,0),U^-)=0$. Similarly, If i.2) holds then $Z_1^-=\{0\}$ and $Z_p^-=\emptyset$ and we can take $U_p^-=\emptyset$ and $U^-=U_0$ in \reftheo{indgen}. From \mref{indf1}, $\mbox{ind}(F_1(\cdot,0),U^-)=\mbox{ind}(F_1(\cdot,0),U_0)=0$.

Hence, $\mbox{ind}(F_1(\cdot,0),U^-)=0$ in both cases.
By \reftheo{indexthm1}, we find a neighborhood $V$ of $0$ in $\mP_2$ such that $\mbox{ind}(T,U\oplus V)=\mbox{ind}(F_1(\cdot,0),U^-)=0$. Since $\mbox{ind}(T,\mB\cap\mX)=1$, we see that $T$ has a fixed point in $\mB\setminus\overline{U\oplus V}$. This fixed point is nontrivial because we are assuming that $T$ has no semi trivial fixed point in $\mP_2$. \eproof

\subsection{Notes on a more special case and a different way to define $T$:}\label{ssapp1}

In many applications, it is reasonable to assume that the cross diffusion effects by other components should be proportional  to the density of a given component. This is to say that if $A(u)=(a_{ij}(u))$ then there are smooth functions $b_{ij}$ such that \beqno{sktcd}a_{ij}(u) = u_ib_{ij}(u) \quad \mbox{if $j\ne i$}.\eeq

In this case, instead of using \mref{b2z}, we can define the $i$-th component $w_i$ of $T(u)$  by
\beqno{Tnew}L_i(u)w_i + \sum_j k_{ij}w_j = f_i(u)+\sum_j k_{ij}u_j,\eeq
where\beqno{L-op}L_i(u)w =-\Div(a_{ii}(u)Dw + w(\sum_{j\ne i}b_{ij}(u)Du_j )).\eeq

As $u\in\mB\cap \mP$, \mref{Tnew} is a {\em weakly coupled} system with H\"older continuous coefficients. We will see that the condition K.1) on the positivity of solutions in the previous section is verified.
To this end, we recall the maximum principles for cooperative linear systems in \cite{FM,FT} and give here an alternative and simple proof to \cite[Theorem 1.1]{FM}. In fact, we consider a more general setting that covers both Dirichlet and Neumann boundary conditions. 

Let us define \beqno{L-op1}L_iw =-\Div(\ag_{i}(x)Dw +\bg_i(x)w) ,\eeq where $\ag_i\in L^\infty(\Og)$, $\bg_i\in L^\infty(\Og,\RR^n)$. Denote $\ccF = (\ccF_1, \cdots, \ccF_m)$.  We then have the following weak minimum principle.

\blemm{FMmax} Let $w$ be a weak solution to the system 
$$\left\{\barr{l}L_i w_i+Kw=\ccF_i,\quad i=1,\ldots,m,  \mbox{ in $\Og$}, \\ \mbox{Homogeneous Dirichlet or Neumann boundary conditions on $\partial \Og$}.\earr\right.$$ 
Assume that $\ag_i(x)\ge\llg_i$ for some $\llg_i>0$ and  $\ccF_i\ge0$ for all $i$. If  $k_{ij}\le0$ for $i\ne j$ and $k_{ii}$ are sufficiently large, in terms of $\sup_\Og\bg_i(u(x))$, then $w\ge0$. \elemm

\bproof Let $\fg^+,\fg^-$ denote the positive and negative parts of a scalar function $\fg$, i.e. $\fg=\fg^+-\fg^-$. We note that $\myprod{Dw_i,Dw_i^-}=-|Dw_i^-|^2$ and $w_iDw_i^-=-w_i^-Dw_i^-$. Integrating by parts, we have
$$\iidx{\Og}{L_iw_iw_i^-} =\iidx{\Og}{(-\ag_i|Dw_i^-|^2-\myprod{\bg_i,Dw_i^-}w_i^-)}$$

Hence, multiplying the $i$-th equation of the system by $-w_i^-$, we obtain $$\sum_i\iidx{\Og}{(\ag_i|Dw_i^-|^2+\myprod{\bg_i,Dw_i^-}w_i^-)}- \sum_{i,j}\iidx{\Og}{k_{ij}w_iw_j^-}=-\sum_i\iidx{\Og}{\myprod{\ccF_i,w_i^-}}. $$

Since $\ccF_i,w_i^-\ge0$, we get $$\sum_i\iidx{\Og}{(\ag_i|Dw_i^-|^2+\myprod{\bg_i,Dw_i^-}w_i^-)}- \sum_{i,j}\iidx{\Og}{k_{ij}w_iw_j^-}\le0. $$

Since $w_iw_i^-=-w_i^-w_i^-$ and   $w_iw_j^-=w_i^+w_j^--w_i^-w_j^-\ge w_i^-w_j^-$,  the above yields, using the assumption that $k_{ij}\le0$ for $i\ne j$
\beqno{kij}\sum_i\iidx{\Og}{(\ag_i|Dw_i^-|^2+\myprod{\bg_i,Dw_i^-}w_i^-+k_{ii}|w_i^-|^2)}- \sum_{i\ne j}\iidx{\Og}{k_{ij}w_i^-w_j^-}\le0.\eeq

By Young's inequality, for any $\eg>0$ we can find a constant $C(\eg,\bg_i)$, depending on $\sup_\Og\bg_i(u(x))$, such that $$\left|\iidx{\Og}{\myprod{\bg_i,Dw_i^-}w_i^-}\right|\le \eg \iidx{\Og}{|Dw_i^-|^2}+C(\eg,\bg_i)\iidx{\Og}{|w_i^-|^2}.$$

Thus, \mref{kij} implies
$$\sum_i\iidx{\Og}{(\ag_i-\eg)|Dw_i^-|^2} +\sum_i (k_{ii}-C(\eg,\bg_i))\iidx{\Og}{|w_i^-|^2}-\sum_{i\ne j}\iidx{\Og}{k_{ij}w_i^-w_j^-}\le0. $$

Combining the ellipticity assumption and Poincar\'e's inequality, we have
$$c_0(\llg_i-\eg)\iidx{\Og}{|w_i^-|^2}\le\iidx{\Og}{(\ag_i-\eg)|Dw_i^-|^2}$$
for some $c_0>0$. Therefore,
$$c_0(\llg_i-\eg)\sum_i\iidx{\Og}{|w_i^-|^2}+\sum_i (k_{ii}-C(\eg,\bg_i))\iidx{\Og}{|w_i^-|^2}-\sum_{i\ne j}\iidx{\Og}{k_{ij}w_i^-w_j^-}\le0. $$

This implies \beqno{cgpos}\iidx{\Og}{\sum_{i,j}\cg_{ij}w_i^-w_j^-}\le0,\eeq where
$$\cg_{ij}=\left\{\barr{ll} c_0(\llg_i-\eg)+k_{ii}-C(\eg,\bg_i)& i=j,\\ -k_{ij} & i\ne j.\earr\right.$$

It is clear that if $k_{ii}$ is sufficiently large then the matrix $\cg=(\cg_{ij})$ is positive definite, i.e. $\myprod{\cg x,x}\ge c|x|^2$ for some positive $c$. Thus, \mref{cgpos} forces $w_i^{-}=0$ a.e. and $w\ge0$. \eproof

Thanks to this lemma, we now see how to construct a matrix $K$ such that $T$ maps $\mB\cap\mP$ into $\mP$. To this end, we note that$f^{(u)}(0,0)=0$ and $ f^{(v)}(u,0)=0$ so that we can write 
$$f^{(u)}(u,v)+K_1^{(u)}u+K_1^{(v)}v=\int_0^1(f_u^{(u)}(tu,tv)+K_1^{(u)})\,dt\,u+\int_0^1(f_v^{(u)}(tu,tv)+K_1^{(v)})\,dt\,v,$$
$$f^{(v)}(u,v)+K_2^{(v)}v= \int_0^1(f_v^{(v)}(u,tv)+K_2^{(v)})\,dt\,v.$$ 

Since $\|\hat{f}_u(u)\|_\mX$ is bounded for $u\in\mB=B(0,M)$ (the bound $M$ is independent of $K$), it is not difficult to see that if the reaction is 'cooperative', i.e., $\partial_{u_i}\hat{f}_j\ge0$ for $i\ne j$, then we can always find $K$ with $k_{ij}=0$ for $i\ne j$ and $k_{ii}>0$ sufficiently large such that the matrix integrands in the above equations are all positive. Therefore,
$\ccF:=\hat{f}(u,v)+K(u,v)\ge0$ for  $(u,v)\in\mB\cap\mP$ and the lemma can apply here.

Finally, for future use in the next section we now explicitly describe the map $T'(u)$ in this case.  For $\Fg=(\fg_1,\ldots,\fg_m)$,  by \reflemm{linlemm}, the components $w_i=T_i(u)$ of $T(u)$ 
and  $W_i=T_i'(u)\Fg$ of $T'(u)\Fg$ solves

\beqno{gen-eig0}
-\Div(\ccA_i(u,W_i)+\ccB_i(u,w_i,\Fg)) +\sum_jk_{ij}W_j=\sum_j(\partial_{u_j}f_i(u)+k_{ij})\fg_j,\eeq where \beqno{gen-eig1}\ccA_i(u,W_i):=a_{ii}(u)DW_i +(\sum_{j\ne i}b_{ij}(u)Du_j )W_i,\eeq
\beqno{gen-eig2}\ccB_i(u,w_i,\Fg):=\sum_j\partial_{u_j}a_{ii}(u)\fg_j Dw_i+w_i\sum_{j\ne i}[\partial_{u_k}b_{ij}(u)\fg_kDu_j +b_{ij}(u)D\fg_j].\eeq

Consider a semi trivial solution $u\in Z_1$, i.e.  for some integer $m_1\ge0$  $$T(u,0)=(u,0),\quad (u,0)=(u_1,\ldots,u_{m_1},0,\ldots,0).$$ 

For $i>m_1$ we have that $w_i=T_i(u,0)$ and $Dw_i=D(T_i(u,0))$ are zero so that $W_i$ solves
\beqno{keyevec}\barr{ll}\lefteqn{-\Div(a_{ii}(u,0)DW_i+\sum_{j\le m_1}b_{ij}(u,0)Du_j W_i)+\sum_j k_{ij}W_j=}\hspace{7cm}&\\& \sum_j(\partial_{u_j}f_i(u,0)+k_{ij})\fg_j.\earr\eeq

\section{Nonconstant and Nontrivial Solutions}\label{nontrivial}\eqnoset

We devote this section to the study of \mref{b1} with Neumann boundary conditions. \reftheo{indgenexist} gives the existence of positive nontrivial solution but this solution may be a constant solution. This is the case when there is a constant vector $u^*=(u^*_1,\ldots,u^*_m)$ such that $f(u^*)=0$. Obviously $u=u_*$ is a nontrivial solution to \mref{b1} and \reftheo{indgenexist} then yields no useful information.
In applications, we are interested in finding a nonconstant solution besides this obvious solution. We will assume throughout this section that the semi trivial solutions are all constant and show that cross diffusion will play an important role for the existence of nonconstant and nontrivial solutions.

Inspired by the SKT systems, we assume that the diffusion is given by \mref{sktcd} as in \refsec{ssapp1} and the reaction term in the $i$-th equation is also proportional to the density $u_i$. This means, \beqno{figi}f_i(u)=u_ig_i(u)\eeq for some $C^1$ functions $g_i$'s. A constant solution $u^*$ exists if it is a solution to the equations $g_i(u^*)=0$ for all $i$.

Throughout this section,  we denote by $\psi_i$'s the eigenfunctions of $-\Dg$, satisfying Neumann boundary condition, to the eigenvalue $\hat{\llg}_i$ such that $\{\psi_i\}$ is a basis for $W^{1,2}(\Og)$. That is,
$$ \left\{\barr{l} -\Delta\psi=\hat{\llg}_i\psi \mbox{ in $\Og$,}\\ \mbox{$\psi$ satisfies homogeneous Neumann boundary condition.}\earr\right.$$

\subsection{Semi trivial constant solutions}\label{sssemiconst} We consider a semi trivial solution $(u,0)$ with $u=(u_1,\ldots,u_{m_1})$ for some integer $m_1=0,\ldots,m$.
Following the analysis of \refsec{ssindex},  we need to consider the eigenvalue problem  $\partial_vF_2(u,0)\Fg_2=\mu\Fg_2$ with $\Fg_2=(\fg_{m_1},\ldots,\fg_m)$ and  $\Fg=(0,\Fg_2)$. Then the equation \mref{keyevec}, with $W_i = \mu\fg_i$ for $i=m_1+1,\ldots,m$, gives
$$-\Div(a_{ii}(u,0)D\fg_i+\sum_{j\le m_1}b_{ij}(u,0)Du_j \fg_i)+\sum_{j>m_1}k_{ij}\fg_j =\mu^{-1}\sum_{j>m_1}(\partial_{u_j}f_i(u,0)+k_{ij})\fg_j.$$

If $u$ is a constant vector then $Du_j=0$ and the above reduces to
\beqno{h1zzzzzz}-\Div(a_{ii}(u,0)D\fg_i)+\sum_{j>m_1}k_{ij}\fg_j =\mu^{-1}\sum_{j>m_1}(\partial_{u_j}f_i(u,0)+k_{ij})\fg_j,\eeq which is an elliptic system with constant coefficients. We then need the following lemma.

\newc{\mbc}{{\mathbf{c}}}

\blemm{eigenlem} Let $A,B$ be constant matrices. Then the solution space of the problem $$\left\{\barr{l}-\Div(AD\Fg) = B\Fg,\\\mbox{Neumann boundary conditions}.\earr\right.$$ has a basis $\{\mbc_{i,j}\psi_i\}$ where  $\mbc_{i,j}$'s are the basis vectors of $\mbox{Ker}(\hat{\llg}_iA-B)$.
\elemm 

The proof of this lemma is elementary. If $\Fg$ solves its equation of the lemma then we can write $\Fg=\sum k_i\psi_i$, in $W^{1,2}(\Og)$, with $k_i\in\RR^m$. We then have $\sum_i\hat{\llg}_iA k_i\psi_i=\sum_iB k_i\psi_i$. Since $\{\psi_i\}$ is a basis of $W^{1,2}(\Og)$, this equation implies $\hat{\llg}_iA k_i=Bk_i$ for all $i$. Thus, $k_i$ is a linear combination of $\mbc_{i,j}$'s.
It is easy to see that $\{\mbc_{i,j}\psi_j\}$ is linearly independent if $\{\mbc_{i,j}\},\{\psi_j\}$ are. The lemma then follows.

\blemm{evectorlem1} Let $m_1$ be a nonnegative integer less than $m$ and $u=(u_1,\ldots,u_{m_1})$ be a constant function such that $T(u,0)=(u,0)$.
Then
$$\partial_{v}f^{(v)}(u,0)\mbc=\llg K_2^{(v)}\mbc$$
has a positive eigenvector $\mbc$ to a positive (respectively, negative) eigenvalue $\llg$ if and only if the eigenvalue problem  $\partial_vF_2(u,0)\Fg_2=\mu\Fg_2$ has a positive solution for some $\mu>1$ (respectively, $\mu<1$).
\elemm

\bproof By \mref{h1zzzzzz} , the eigenvalue problem  $\partial_vF_2(u,0)\Fg_2=\mu\Fg_2$ (or $W_i=\mu\fg_i$) is determined by the following system
\beqno{h1zzzzzz1}-\mu\Div(A^{(m_1)}(u,0)D\Fg_2)+\mu K_2^{(v)}\Fg_2 =[\partial_{v}f^{(v)}(u,0)+K_2^{(v)}]\Fg_2,\eeq  where   $A^{(m_1)}(u,0)=\mbox{diag}[a_{ii}(u,0)]_{i>m_1}$. The coefficients of the above system are constant and \reflemm{eigenlem} yields that the solutions to the above is 
$\sum \mbc_i\psi_i$ with $\mbc_i$ solving
$$\mu[\hat{\llg}_iA^{(m_1)}(u,0)+K_2^{(v)}]\mbc =[\partial_{v}f^{(v)}(u,0)+K_2^{(v)}]\mbc.$$ 

Note that the only positive eigenfunction of $-\Delta$ is $\psi_0=1$ to the eigenvalue $\hat{\llg}_0=0$. Therefore, from the above system with $i=0$ we see that if the system $$\mu K_2^{(v)}\mbc =[\partial_{v}f^{(v)}(u,0)+K_2^{(v)}]\mbc \Leftrightarrow \partial_{v}f^{(v)}(u,0)\mbc=(\mu-1)K_2^{(v)}\mbc$$ has a positive solution $\mbc$ then the constant function $\mbc$ is a positive eigenfunction for $\partial_vF_2(u,0)$. Conversely, if $\partial_vF_2(u,0)\Fg_2=\mu\Fg_2$ has a positive solution $\Fg_2$ then we integrate \mref{h1zzzzzz1} over $\Og$, using the Neumann boundary conditions, to see that $\mbc=\iidx{\Og}{\Fg_2}$ is a positive solution to the above system.  The lemma then follows. \eproof

\brem{LVrem1} By the Krein-Ruthman theorem, if $\partial_vF_2(u,0)$ is strongly positive then $\mu=r_v(u)$ is the only eigenvalue with positive eigenvector. The eigenvalue problem  $\partial_vF_2(u,0)\Fg_2=\mu\Fg_2$ has a positive solution for $\mu=1$  if and only if the matrix $\partial_{v}f^{(v)}(u,0)$ has a positive eigenvector to the zero eigenvalue.\erem

We now discuss the special case $f_i(u)=u_ig_i(u)$.

\blemm{LVlem1} Assume that $f_i(u)=u_ig_i(u)$. Let $m_1$ be a nonnegative integer less than $m$ and $u=(u_1,\ldots,u_{m_1})$ be a constant vector such that $T(u,0)=(u,0)$.
Then the eigenvalue problem $\partial_vF_2(u,0)\Fg_2=\mu\Fg_2$ 
\bdes\item[i.1)] has no nonzero solution for $\mu=1$ if and only if  $g_i(u,0)\ne 0$ for any $i>m_1$;

\item[i.2)] has a positive solution  for some $\mu>1$ if and only if   $g_i(u,0)>0$ for some $i>m_1$; 
\item[i.3)] has no positive solution for $\mu>1$ if and only  $g_i(u,0)< 0$ for any $i>m_1$.\edes

\elemm

\bproof We now let $K=kI$. By \reflemm{evectorlem1} the existence of positive eigenvectors of $\partial_vF_2(u,0)\Fg_2=\mu\Fg_2$ is equivalent to that of \beqno{fvev}\partial_{v}f^{(v)}(u,0)\mbc=k\llg \mbc \quad \mbox{with $\mbc=(c_{m_1+1},\ldots,c_m)>0$ and $\llg=\mu-1$}.\eeq

Since $f_i(u)=u_ig_i(u)$ and $u_i=0$ for $i>m_1$, we have $\partial_{u_k}f_i(u,0)=\dg_{ik}g_i(u,0)$, where $\dg_{ik}$ is the Kronecker symbol, for $i,k>m_1$. Thus, $\partial_{v}f^{(v)}(u,0)$ is a diagonal matrix and \mref{fvev} is simply $$g_i(u,0)c_i=k\llg c_i \quad \forall i>m_1.$$ 

Clearly i.1) holds because then the above system has nonzero eigenvector to $\llg=0$. 
For i.2)  we can take $\llg=g_i(u,0)/k>0$ and $c_i=1$, other components of $\mbc$ can be zero.  i.3) is obvious. The proof is complete. \eproof

We then have the following theorem for systems of two equations.
\btheo{LV2thm1} Assume that $f_i(u)=u_ig_i(u)$ for $i\in\{1,2\}$. Suppose that the trivial and semi trivial solutions are only the {\em constant} ones $(0,0)$, $u_{1,*}$
and $u_{2,*}$. This means, $g_i(u_{i,*})=0$.Then there is a nontrivial solution $(u_1,u_2)>0$ in the following situations:
\bdes\item[a)] $g_i(0)> 0$, $i=1,2$, and $g_1(u_{2,*})$ and $g_2(u_{1,*})$ are positive.
\item[b)] $g_i(0)> 0$, $i=1,2$, and $g_1(u_{2,*})$ and $g_2(u_{1,*})$ are negative.
\item[c)] $g_1(0)> 0$, $g_2(0)<0$, and $g_2(u_{1,*})>0$.
\edes \etheo

\bproof We just need to compute the local indices of $T$ at the trivial and semi trivial solutions and show that the sum of these indices is not 1. 

First of all, by i.1) of \reflemm{LVlem1}, it is clear that the condition Z) at these solutions are satisfied in the above situations.

The conditions in case a) and i.2) of \reflemm{LVlem1} imply that 0 and the semi trivial solutions are unstable in theirs complement directions. \reftheo{indexthm1}, with $Z_1^-=\emptyset$, gives that the local indices at these solutions are all zero. Similarly, in case b),
the local index at 0 is 0 and the local indices at the semi trivial solutions, which are stable  in theirs complement directions, are 1. In these cases, the sum of the indices is either 0 or 2.

In case c), because $g_2(0)< 0$ we see that 0 is $u_2$-stable so that $T_2:=T|_{\mX_2}$, the restriction of the map $T$ to $\mX_2$, has its local  index at 0 equal 1 and therefore its local index at $u_{2,*}$ is zero (see also the proof of \reftheo{indgenexist}). The assumption $g_1(0)> 0$ also yields a neighborhood $V_1$ in $\mX_1$ of 0 such that $\mbox{ind}(T,V_1\oplus \mX_2)=0$ (the stability of $u_{2,*}$ in the $u_1$ direction does not matter). On the other hand, because $g_1(0)> 0$ we see that 0 is $u_1$-unstable so that $T_1:=T|_{\mX_1}$, the restriction of the map $T$ to $\mX_1$, has its local index at 0 equal 0 and therefore its local index at $u_{1,*}$ is 1. But $u_{1,*}$ is $u_2$-unstable, because $g_2(u_{1,*})>0$,  so that there is a neighborhood $V_2$ in $\mX_2$ of 0 such that $\mbox{ind}(T,\mX_1\oplus V_2)=0$.

In three cases, we have shown that the sum of the local indices at the trivial and semi trivial solutions  is not 1.  Hence, there is a positive nontrivial fixed point $(u_1,u_2)$.
\eproof

\brem{nontrivrem1} If the system $g_i(u)=0$, $i=1,2$, has no positive constant solution then the above theorem gives conditions for the existence of nonconstant and nontrivial solutions. This means pattern formations occur. \erem

\subsection{Nontrivial constant solutions}\label{ssnotrivconst}
Suppose now that $u_*=(u_1,\ldots,u_m)$ is a nontrivial constant fixed point of $T$ with $u_i\ne0$ for all $i$. We will use the Leray Schauder theorem to compute the local index of $T$ at $u_*$. Since $u_*$ is in the interior of $\mP$, we do not need that $T$ is positive as in the previous discussion so that we can take $K=0$. The main result of this section, \reftheo{sumdim}, yields a formula to compute the indices at nontrivial constant fixed points. In applications, the sum of these indices and those at semi trivial fixed points will shows the possibility of nontrivial and nonconstant fixed points to exist. 

In the sequel, we will denote \beqno{diagA}d_A(u_*)=\mbox{diag}[a_{11}(u_*),\ldots,a_{mm}(u_*)].\eeq From the ellipticity assumption on $A$, we easily see that $a_{ii}(u_*)>0$ for all $i$ and thus $d_A(u_*)$ is invertible.

The following lemma describes the eigenspaces of $T'(u_*)$.

\blemm{evectorlem2}  The solution space of $T'(u_*)\Fg=\mu\Fg$ is spanned by $\mbc_i\psi_i$ with $\mbc_i$ solving
\beqno{ev1}\hat{\llg}_i[A(u_*)+(\mu-1)d_A(u_*)]\mbc_i=\partial_u F(u_*)\mbc_i.\eeq 

\elemm

\bproof We have $D(T(u_*))=Du_*=0$ so that \mref{gen-eig0}-\mref{gen-eig2}, with $w=u_*$, show that the $i$-th component $W_i$ of $T'(u_*)\Fg$ solves 
$$-\Div(a_{ii}(u_*)DW_i+w_i\sum_{j\ne i}b_{ij}(u_*)D\fg_j ) =\partial_{u_k}f_i(u_*)\fg_k.$$ 

Since $w_ib_{ij}(u_*) =u_ib_{ij}(u_*)= a_{ij}(u_*)$,  the eigenvalue problem $W=T'(u_*)\Fg=\mu\Fg$ is
$$-\Div(a_{ii}(u_*)D(\mu\fg_i)+\sum_{j\ne i}a_{ij}(u_*)D\fg_j ) =\partial_{u_k}f_i(u_*)\fg_k.$$

In matrix form,  the above can be written as \beqno{ev0}-\Div([A(u_*)+(\mu-1)d_A(u_*)]D\Fg) = \partial_u F(u_*)\Fg.\eeq

Since $u_*$ is a constant vector, by \reflemm{eigenlem} we can write $\Fg=\sum \mbc_i\psi_i$ with $\mbc_i$ solving $$\hat{\llg}_i[A(u_*)+(\mu-1)d_A(u_*)]\mbc_i=\partial_u F(u_*)\mbc_i \quad \forall i.$$ This is \mref{ev1} and the lemma is proved. \eproof

We now have the following  explicit formula for $\mbox{ind}(T,u_*)$.
\btheo{sumdim} Assume that \beqno{No1}\mbox{Ker}(\hat{\llg}_i A(u_*)-\partial_u F(u_*))=\{0\} \quad \forall i.\eeq
For $i>0$ let   $\ccA_i=A(u_*)-\hat{\llg}_i^{-1}\partial_u F(u_*)$ and  $$N_{i}=\sum_{\llg<0}\mbox{dim}(\mbox{Ker}(d_A(u_*)^{-1}\ccA_i-\llg I)).$$ We also denote by $M_i$ the multiplicity of $\hat{\llg}_i$.

Then there exists an integer $L_0$ such that \beqno{dimdef} \mbox{ind}(T,u_*)=(-1)^\cg, \quad\cg=\sum_{i\le L_0} N_{i}M_i.\eeq  \etheo

\bproof We will apply Leray-Schauder's theorem to compute $\mbox{ind}(T,u_*)$.  First of all,  \reflemm{evectorlem2} and \mref{No1}  show that $\Fg=0$ is the only solution to $T'(u_*)\Fg=\Fg$  so that $\mu=1$ is not an eigenvalue of $T'(u_*)$.

By Leray-Schauder's theorem, we have that $\mbox{ind}(T,u_*)=(-1)^\cg$, where $\cg$ is the sum of multiplicities of eigenvalues $\mu$ of $T'(u_*)$ which are greater than 1. \reflemm{evectorlem2} then clearly shows that $\cg$ is the sum of the dimensions of solution spaces of \mref{ev1} and  \beqno{mult1}\cg:=\sum_{i} \cg^*_{i}M_i,\eeq  and $\cg^*_i=\sum_{\mu>1} n_{i,\mu}$, where $n_{i,\mu}$ is the dimension of the solution space of $$\hat{\llg}_i[A(u_*)+(\mu-1)d_A(u_*)]\mbc-\partial_u F(u_*)\mbc=0.$$

For $i=0$ we have $\hat{\llg}_0=0$ so that $n_{0,\mu}=\mbox{dim}(\mbox{Ker}(\partial_u F(u_*)))$, which is zero because of \mref{No1}.
For $i>0$ and $\mu>1$ let   $\ccA_i=A(u_*)-\hat{\llg}_i^{-1}\partial_u F(u_*)$, as being defined in this theorem, and $\llg=1-\mu<0$. The above equation can be rewritten as $ \ccA_i\mbc=\llg d_A(u_*)\mbc $ so that
$\mbc$ is an eigenvector to a negative eigenvalue $\llg$ of the matrix $d_A(u_*)^{-1}\ccA_i$.
It is clear that the number $\cg^*_i$ in \mref{mult1} is the sum of the dimensions of eigenspaces of $d_A(u_*)^{-1}\ccA_i$ to negative eigenvalues. That is $\cg^*_i=N_i$.

From the ellipticity assumption on $A$, it is not difficult to see that $d_A(u_*)^{-1}A(u_*)$ is positive definite. Therefore,  $d_A(u_*)^{-1}\ccA_i$ is positive if $\hat{\llg}_i$ is large. Thus, as $\lim_{i\to\infty}\hat{\llg}_i=\infty$, there is an integer $L_0$ such that $d_A(u_*)^{-1}\ccA_i$ has no negative eigenvalues  if $i> L_0$.  Hence, $\cg^*_{i}=N_i=0$ if $i> L_0$.  The theorem then follows.
\eproof

\brem{No1rem} Since $\hat{\llg}_0=0$, \mref{No1} implies that $\mbox{Ker}(\partial_u F(u_*))=\{0\}$ so that $u_*$ is an isolated constant solution to $F(u_*)=0$.  Also, as $A(u_*)$ is positive definite and $\lim_{i\to\infty}\hat{\llg}_i=\infty$, we see that \mref{No1} is true when $i$ is large. \erem

Combining \reftheo{sumdim} and \reftheo{LV2thm1}, we obtain
\bcoro{LV2coro1} Assume as in \reftheo{LV2thm1}. Suppose further that there is only one nontrivial constant solution $u_*$. Let $\cg$ be as in \mref{dimdef}. There is a nontrivial nonconstant solution $(u_1,u_2)>0$ in the following situations:
\bdes\item[1)] $\cg$ is odd and a) or c) of the theorem hold.
\item[2)] $\cg$ is even and b) of the theorem holds.
\edes \ecoro

\bproof We have seen from the proof of \reftheo{LV2thm1} that the sum of the local indices at the trivial and semi trivial solutions is 0 in the cases a) and c). If $\cg$ is odd then $\mbox{ind}(T,u_*)=-1$. Similarly, if b) holds then the sum of the indices at the trivial and semi trivial solutions is 2. If $\cg$ is even then $\mbox{ind}(T,u_*)=1$. Thus,  the sum of the local indices at constant solutions is not 1 in both cases. Since $\mbox{ind}(T,\mB\cap\mX)=1$, a nonconstant and  nontrivial solution must exist. \eproof

\subsection{Some nonexistent results}\label{ssnonext}
We conclude this paper by some nonexistence results showing that if the parameter $\llg_0$ is sufficiently large then there is no nonconstant solutions. We consider the following system
\beqno{b1aaa}\left\{\barr{ll} -\Div(A(u,Du)=f(u)+B(u,Du)& \mbox{in $\Og$},\\
\mbox{homogenenous Neumann boundary conditions}& \mbox{on $\partial \Og$}, \earr\right.\eeq 

We first have the following nonexistent result under a strong assumption on the uniform boundedness of solutions. This assumption will be relaxed later in \refcoro{advnonext}.
\btheo{nonextthm} Assume that $A$ satisfies A) and $\hat{f}(u,Du):=f(u)+B(u,Du)$ for some $f\in C^1(\RR^m,\RR^m)$ and $B:\RR^m\times\RR^{mn}\to\RR^m$ such that $|B(u,p)|\le b(u)|p|$ for some continuous nonnegative function $b$ on $\RR^m$. Suppose also that there is a constant $C$ independent of $\llg_0$ such that for any solutions of \mref{b1aaa}  $$\|u\|_{L^\infty(\Og)}\le C.$$ If the constant $\llg_0$ in A) is sufficiently large then there is no nonconstant solution to \mref{b1aaa}. \etheo

\bproof For any function $g$ on $\Og$ let us denote  the average of $g$ over $\Og$ by $g_\Og$. That is,
$$g_\Og =\frac{1}{|\Og|}\iidx{\Og}{g}.$$

Integrating \mref{b1aaa} and using Neumann boundary conditions, we have $f(u)_\Og +B(u,Du)_\Og=0$. Thanks to this, we test the system with $u-u_\Og$ to get \beqno{nonc1}\iidx{\Og}{\myprod{A(u,Du),Du}} = \iidx{\Og}{[\myprod{f(u)-f(u)_\Og,u-u_\Og}-\myprod{B(u,Du)_\Og,u-u_\Og}]}.\eeq

We estimate the terms on the right hand side. First of all, by
H\"older's inequality 
$$\iidx{\Og}{\myprod{f(u)-f(u)_\Og,u-u_\Og}} \le \left(\iidx{\Og}{|f(u)-f(u)_\Og|^2}\right)^\frac12
\left(\iidx{\Og}{|u-u_\Og|^2}\right)^\frac12.$$
Applying Poincar\'e's inequality to the functions $f(u),u$ on the right hand side of the above inequality, we can bound it by
$$ \mbox{diam}(\Og)^2\left(\iidx{\Og}{|f_u(u)|^2|Du|^2}\right)^\frac12
\left(\iidx{\Og}{|Du|^2}\right)^\frac12\le F_*\mbox{diam}(\Og)^2\|Du\|_{L^2(\Og)}^2,$$ where $F_*:=\sup_{\Og}|f_u(u(x))|$. This number is finite because we are assuming that $\|u\|_{L^\infty(\Og)}$ is bounded uniformly.

Similarly, we define $B_*:=\sup_{\Og}b(u(x))$. Using the facts that $|B(u,Du)|\le b(u)|Du|\le B_*|Du|$, we have $|B(u,Du)_\Og|\le B_*|\Og|^{-1}\|Du\|_{L^1(\Og)}$. Furthermore, by H\"older's and Poincar\'e's inequalities
$$\|Du\|_{L^1(\Og)}\le |\Og|^\frac12 \|Du\|_{L^2(\Og)},\; \|u-u_\Og\|_{L^1(\Og)}\le C|\Og|^\frac12 \mbox{diam}(\Og) \|Du\|_{L^2(\Og)}.$$ We then obtain $$\left|\iidx{\Og}{\myprod{B(u,Du)_\Og,u-u_\Og}}\right|\le B_*\mbox{diam}(\Og) \|Du\|_{L^2(\Og)}^2.$$

Using the above estimates and the ellipticity condition A) in \mref{nonc1}, we get 
$$\llg_0\iidx{\Og}{|Du|^2}\le F_*\mbox{diam}(\Og)^2
\iidx{\Og}{|Du|^2}+B_*\mbox{diam}(\Og)
\iidx{\Og}{|Du|^2}.$$

If $\llg_0$ is sufficiently large then the above inequality clearly shows that $\|Du\|_{L^2(\Og)}=0$ and thus $u$ must be a constant vector. \eproof

\brem{DBCnoex} If we assume Dirichlet boundary conditions and $f(0)\equiv0$ then 0 is the only solution if $\llg_0$ is sufficiently large. To see this we test the system with $u$ and repeat the argument in the proof. \erem

The assumption on the boundedness of the $L^\infty$ norms of the solutions in \reftheo{nonextthm} can be weakened if $\llg(u)$ has a polynomial growth. We have the following result. 

\bcoro{advnonext}  The conclusion of \reftheo{nonextthm} continues to hold for the system
if  one has a uniform estimate for $\|u\|_{W^{1,2}(\Og)}$ and $\llg(u)\sim \llg_0 + (1+|u|)^k$ for some $k>0$. \ecoro

\bproof We just need to show that the two assumptions in fact provide the uniform bound of $L^\infty$ norms needed in the previous proof. From the growth assumption on $\llg$, we see that the number $\mathbf{\LLg}$ in \mref{LLg} of \refsec{estsec} is now $$\mathbf{\LLg}=\sup_{W\in\RR^m}\frac{|\llg_W(W)|}{\llg(W)}\sim \sup_{W\in\RR^m}\frac{(1+|W|)^{k-1}}{\llg_0+(1+|W|)^{k}}.$$ By considering the cases $(1+|W|)^{k}$ is greater or less that $\llg_0$, we can easily see that $\mathbf{\LLg}$ can be arbitrarily small if $\llg_0$ is sufficiently large. On the other hand, our assumptions yield that $\|u\|_{L^1(\Og)}$ is uniformly bounded. Thus, we can fix a $R_0>0$ and use the fact that $\|u\|_{BMO(B_{R_0})}\le C(R_0)\|u\|_{L^1(\Og)}$ to see that the condition D) of \refprop{Dup-uni} holds if $\llg_0$ is sufficiently large. We then have that the H\"older norms, and then $L^\infty$ norms, of the solutions to \mref{b1aaa} are uniformly bounded, independently of $\llg_0$. This is the key assumption of the proof of \reftheo{nonextthm} so that the proof can continue as before. \eproof

\brem{Xuniest1} From the examples in \refsec{ssstrong}, e.g. \reflemm{Ducontrollem} and \refcoro{nis2}, we see that the assumption on uniform bound for $W^{1,2}(\Og)$ norms can be further weakened by the same assumption on $L^1(\Og)$ norms.  \erem

\bibliographystyle{plain}

\end{document}